\documentclass[12pt,oneside]{amsart}
\usepackage[T1]{fontenc}
\usepackage{lmodern}
\input{glyphtounicode}
\usepackage{amsmath,amssymb,amsthm,mathrsfs}
\usepackage[margin=1.05in]{geometry}
\usepackage{enumitem}
\usepackage{needspace}
\usepackage[hidelinks]{hyperref}
\newtheorem{thm}{Theorem}[section]
\newtheorem{prop}[thm]{Proposition}
\newtheorem{lem}[thm]{Lemma}
\newtheorem{cor}[thm]{Corollary}
\theoremstyle{definition}
\newtheorem{defn}[thm]{Definition}
\newtheorem{example}[thm]{Example}
\theoremstyle{remark}
\newtheorem{rem}[thm]{Remark}
\numberwithin{equation}{section}
\newcommand{\norm}[1]{\left\|#1\right\|}
\newcommand{\abs}[1]{\left|#1\right|}
\newcommand{\eps}{\varepsilon}
\newcommand{\C}{\mathbb C}

\newcommand{\N}{\mathbb N}
\newcommand{\T}{\mathrm{T}}
\newcommand{\Aff}{\operatorname{Aff}}

\newcommand{\Tr}{\operatorname{Tr}}
\newcommand{\id}{\operatorname{id}}

\newcommand{\calF}{\mathcal F}
\newcommand{\calD}{\mathcal D}
\newcommand{\calM}{\mathcal M}
\newcommand{\Prob}{\operatorname{Prob}}

\newcommand{\B}{\mathcal B}
\newcommand{\Aprob}{\mathfrak A}
\newcommand{\Probfa}{\operatorname{Prob}_{\mathrm{fa}}}
\newcommand{\Proj}{\operatorname{Proj}}
\newcommand{\Stone}{\operatorname{Stone}}

\newcommand{\symdiff}{\mathbin{\triangle}}
\setlist[enumerate]{label=\textup{(\roman*)},leftmargin=*}
\setlist[itemize]{leftmargin=*}
\hypersetup{
  pdftitle={Classification of uniform tracial Cartan pairs for amenable group actions},
  pdfauthor={George A. Elliott and Zhuang Niu},
  pdfkeywords={Cartan subalgebras, uniform tracial completion, small boundary property, orbit equivalence, amenable group actions}
}
\title[Uniform tracial Cartan pairs]{Classification of uniform tracial Cartan pairs for amenable group actions}
\author{George A. Elliott}
\address{Department of Mathematics, University of Toronto,
40 St. George Street, Toronto, Ontario M5S 2E4, Canada}
\email{elliott@math.toronto.edu}
\author{Zhuang Niu}
\address{Department of Mathematics and Statistics, University of Wyoming,
Dept. 3036, 1000 E. University Avenue, Laramie, Wyoming 82071, USA}
\email{zniu@uwyo.edu}
\subjclass[2020]{Primary 46L05; Secondary 37A20, 37A55, 46L40}
\keywords{Cartan subalgebras, uniform tracial completion, small boundary
property, uniform measure orbit equivalence, amenable group actions}
\date{}
\begin{document}
\begin{abstract}
Consider free and minimal actions of countably infinite amenable groups
acting on compact metrizable spaces with small boundary property.
The uniform tracial completions of the canonical Cartan pairs of the
crossed product $\mathrm{C}^*$-algebras are classified by their trace simplexes.
Moreover, every affine homeomorphism between the trace simplexes is
induced by an orbit equivalence on Borel sets conull for every invariant
probability measure.
\end{abstract}
\maketitle
\raggedbottom

\section{Introduction}

In 1959 Dye showed that any two ergodic measure-preserving automorphisms of
standard nonatomic probability spaces are orbit equivalent modulo null
sets~\cite{Dye}. Ornstein and Weiss extended this conclusion to free
ergodic probability-preserving actions of countably infinite amenable
groups~\cite{OW}. Connes, Feldman, and Weiss showed that every amenable
countable measured equivalence relation is generated, modulo null sets,
by a single transformation. They deduced that any two Cartan subalgebras
of the hyperfinite $\mathrm{II}_1$ factor are conjugate by an
automorphism~\cite{CFW}.

For minimal homeomorphisms of Cantor sets, Giordano, Putnam, and Skau
characterized topological orbit equivalence by the existence of a
homeomorphism of the underlying spaces carrying one family of invariant
probability measures onto the other~\cite[Theorem~2.2]{GPS95}.
Here we ask whether a prescribed affine homeomorphism between the
invariant-measure simplexes can be realized by a single orbit
equivalence on invariant Borel sets conull for every invariant
probability measure. We prove this for free minimal actions of countably
infinite amenable groups on compact metrizable spaces with the small
boundary property.

The same affine homeomorphism is implemented by an isomorphism of the
uniform tracial completions of the canonical Cartan pairs of the crossed
products. In particular, if two such pairs have the same ambient
completion and the same designated traces, the identity on that trace simplex is
implemented by an automorphism carrying one completed diagonal onto
the other and fixing each designated trace. The small boundary
assumption is substantive even in a fixed uniform tracial completion.
Two nonconjugate uniform Cartan subalgebras can occur in a common
completion~\cite[Remark~8.9]{UR}.

Throughout the dynamical parts of the paper, $\Gamma$ is a countably
infinite discrete amenable group acting freely and minimally on a
nonempty compact metrizable space $X$. We use left actions, write $gx$
for the action, and set
\[
 A=\mathrm{C}(X)\rtimes_r\Gamma,\qquad D=\mathrm{C}(X),\qquad
 \Delta=M_\Gamma(X).
\]
The conditional expectation onto $D$ identifies $\Delta$ with $\T(A)$.
The small boundary property means that $X$ has a basis of open sets
$U$ with $\mu(\partial U)=0$ for every $\mu\in\Delta$.
Freeness is used here, as well as in the tower construction. The trace
completion is always taken on operator-norm bounded sets; it is not the
Hilbert-space completion in a $2$-norm.

\Needspace{10\baselineskip}
\begin{thm}\label{thm:main}
Let $\Gamma_i\curvearrowright X_i$, $i=1,2$, satisfy the preceding
conventions and have the small boundary property. Put
$A_i=\mathrm{C}(X_i)\rtimes_r\Gamma_i$ and $D_i=\mathrm{C}(X_i)$. Every affine
homeomorphism $H:\T(A_1)\to\T(A_2)$ is implemented by an isomorphism
\[
 \Theta:(D_1^u,A_1^u)\longrightarrow(D_2^u,A_2^u)
\]
such that $(H\tau)(\Theta(a))=\tau(a)$ for $a\in A_1^u$.
Conversely, a marked tracial isomorphism of these pairs induces an
affine homeomorphism of their trace simplexes.
\end{thm}

\begin{rem}\label{rem:intro-Z-stable}
Theorem~\ref{thm:main} and its orbit-equivalence consequences also apply
when each action has the uniform Rokhlin property \textup{(URP)} and
its crossed product is $\mathcal Z$-stable.
Section~\ref{sec:AF-classification} shows that these hypotheses imply
the small boundary property.
\end{rem}

\begin{cor}\label{cor:intro-uniform-measure-oe}
Under the hypotheses of Theorem~\ref{thm:main}, every affine
homeomorphism
\[
 H:M_{\Gamma_1}(X_1)\to M_{\Gamma_2}(X_2)
\]
is implemented by a
uniform measure orbit equivalence. In particular, there are invariant
universally conull Borel sets $X_i^0$ and a Borel bijection
$\theta:X_1^0\to X_2^0$ with
\[
 \theta(\Gamma_1x)=\Gamma_2\theta(x),\qquad
 \theta_*\mu=H\mu\quad(\mu\in M_{\Gamma_1}(X_1)).
\]
It induces an isometric isomorphism of the uniform measure full groups
by $T\mapsto\theta T\theta^{-1}$.
\end{cor}

The orbit equivalence acts simultaneously over the whole invariant-measure
simplex. Its pullbacks of continuous observables and the indicators of
its cocycle partitions are uniformly approximable in $L^2$ over that
simplex by bounded continuous functions. These are the regularity
requirements in Definition~\ref{defn:uniform-measure-oe}.
Section~\ref{sec:uniform-orbits} proves the corollary.

The following is a straightforward corollary.

\begin{cor}[At most countably many ergodic measures]
\label{cor:countable-ergodic-measures}
Let $\Gamma_i\curvearrowright X_i$, $i=1,2$, be free minimal actions of
countably infinite amenable groups on compact metrizable spaces, each
admitting at most countably many ergodic invariant probability measures.
The actions are uniformly measure orbit equivalent if and only if the
trace simplexes of their crossed products are affinely homeomorphic.
\end{cor}

In particular, among actions with finitely many ergodic invariant
measures, their number determines the uniform measure orbit-equivalence
class.

We also introduce a uniform version of the measured full group.
For $K=M_\Gamma(X)$, the \emph{uniform measure full group}
$[R_\Gamma]^u_K$ consists of orbit-preserving Borel bijections on
invariant sets conull for every $\mu\in K$, modulo equality almost
everywhere for every such measure. We require the indicators of the
displacement sets $\{x:T(x)=gx\}$ to be uniformly approximable in
$L^2(K)$ by continuous functions; see
Definition~\ref{defn:uniform-measure-full-group}. Its complete
bi-invariant metric is
\[
 d_K(T,S)=\sup_{\mu\in K}\mu\{x:T(x)\ne S(x)\}.
\]
Every element preserves every measure in $K$. The canonical orbit
unitaries identify this group with
$\mathcal N_{A^u}(D^u)/\mathcal U(D^u)$, the completed unitary
normalizer modulo diagonal unitaries. Uniform measure orbit
equivalences conjugate these groups isometrically. Conversely, for the
actions considered here, every topological isomorphism between their
uniform full groups is uniquely
spatial and automatically isometric
\cite[Theorem~8.17]{UR}. Thus the topological group determines the
uniform measure orbit-equivalence class
(Remark~\ref{rem:uniform-full-group-reconstruction-action}).
For a single measure, the definition recovers the usual measured full group;
for the whole simplex it imposes simultaneous continuous
approximability. The AF model supplies a dense
locally finite subgroup of this full group, as proved in
Corollary~\ref{cor:uniform-full-group-locally-finite}.
In the Cantor models, finite-order elements of the topological full
group are also dense for $d_K$. These results give a spatial form of
the finite-permutation approximations used in the construction.

There is also a uniform analogue of the single-generator conclusion of
Connes--Feldman--Weiss. Corollary~\ref{cor:uniform-CFW-action} shows that
every system considered here is uniformly measure orbit equivalent to
a Cantor minimal $\mathbb Z$-system. Transporting its generator gives
one element of the uniform measure full group whose orbits are the
original $\Gamma$-orbits on a single invariant Borel set conull for
every invariant probability measure. The equivalence retains both the
completed observables and the cocycle-fiber projections, simultaneously
over the entire invariant-measure simplex.

The norm structure of the completed diagonal is examined in
Section~\ref{sec:boolean-completions}. Bounded uniform completeness
upgrades finite-spectrum approximation to norm real rank zero
(Lemma~\ref{lem:boolean-RR0}), giving the equivalence between
relative SBP and real rank zero of the entire completion
(Theorem~\ref{thm:boolean-SBP-RR0}). This answers
\cite[Remark~2.18]{EN}. The projection Boolean algebra therefore
determines the unmarked norm diagonal. Its exact capacity recovers
the invariant-measure simplex
(Theorem~\ref{thm:boolean-recovery}), and the corresponding
isometric Boolean algebras are classified by that simplex.
The Boolean algebra is $\sigma$-complete exactly when the simplex
has finitely many extreme points
(Theorem~\ref{thm:boolean-completeness}). Consequently, with
infinitely many ergodic invariant measures the completed diagonal
is not $\mathrm{C}^*$-isomorphic to the $L^\infty$-algebra of a
standard probability space. Infinite-vertex simplexes can also
give different bare Boolean algebras, and hence nonisomorphic
norm $\mathrm{C}^*$-algebras
(Proposition~\ref{prop:boolean-different}). We also prove a general
upper-probability range theorem and a Boolean-core criterion for
realization by free minimal Cantor actions. The general range
criterion applies to arbitrary designated compact families;
the simplex-recovery and sharp completeness statements use the
full-trace AF models.

The proof has two parts. First, Theorem~\ref{thm:AF-model} constructs a
simple unital AF algebra with its standard diagonal whose completed
Cartan pair models the dynamical pair. The construction uses the
zero-dimensional extension and clopen tower results of
Kerr--Szab\'o~\cite{KS}, and permits arbitrarily large Bratteli
multiplicities. Second, Theorem~\ref{thm:AF-pairs} matches regular
finite-dimensional maps after deleting corners of uniformly small
trace. The resulting intertwining retains matrix-unit normalizers and
the completed diagonal, and its invariant is the trace simplex.
No Cuntz-comparison assumption is required. The classification is taken
in the bounded uniform tracial topology.

The theory of tracially complete $\mathrm{C}^*$-algebras developed in
\cite{CCEGSTW} provides the ambient setting for this question. The
relative construction here also preserves a specified diagonal and
its normalizers. These data make it possible to realize the
isomorphism spatially as a uniform measure orbit equivalence and to
identify the associated full groups.

The framework for regular uniform measure relations is developed
in~\cite[Section~3]{UR}.
For a compact convex designated family of invariant probabilities,
relative small boundaries and amenability for each designated measure
give compatible finite matrix charts and an AF Cartan model
\cite[Theorems~6.1 and~12.6]{UR}. If the family consists of all invariant
probabilities and the classes are infinite on a common conull set,
the completed pairs and the uniform measure relations are classified
by that simplex \cite[Theorem~8.6]{UR}. For the reconstruction of
regular twisted relation models from completed Cartan pairs and the
description of their isomorphisms, see~\cite[Theorems~15.9 and~15.12]{UR}.
The treatment of controlled quasi-invariant families appears
in~\cite[Section~14]{UR}. Without small boundaries, relation URP is
characterized by uniform local AH Cartan approximation
\cite[Theorem~16.13]{UR}. For actions, the proofs of the Cartan-pair
classification, its orbit-equivalence corollary, and the Boolean
classification appear in Sections~\ref{sec:AF-classification},
\ref{sec:uniform-orbits}, and~\ref{sec:boolean-completions}, respectively.

Section~\ref{sec:preliminaries} collects the uniform tracial
completion conventions and bounded density facts, and establishes
the expectation and Cartan-pair preliminaries.
Section~\ref{sec:cantor-reduction} reduces the dynamics to a Cantor
system. Sections~\ref{sec:AF-model} and~\ref{sec:AF-classification}
construct the AF model and classify the completed AF pairs.
The Borel realization in Section~\ref{sec:uniform-orbits} gives
orbit equivalence, full-group density, and the single-generator
consequence. Section~\ref{sec:interval-maps} proves approximate
conjugacy for diagonal interval maps by unitaries which almost
normalize the diagonal in the uniform $2$-norm.
Section~\ref{sec:boolean-completions} proves norm real rank zero and
develops the Boolean-algebra range, classification, completeness,
and action-realization results.

\medskip
\noindent\textbf{Acknowledgements.}
The work was carried out using ChatGPT.
The first author is supported by an NSERC Discovery Grant, and the
second author is supported by a Simons Foundation Grant (MP-TSM-00002606).

\section{Preliminaries}\label{sec:preliminaries}

\subsection{Uniform tracial completions}\label{subsec:tracial-completions}

\begin{defn}\label{def:completion}
Let $C$ be a unital $\mathrm{C}^*$-algebra and let $K\subseteq\T(C)$ be nonempty,
compact, and convex. Write
\[
 \norm{c}_{2,K}=\sup_{\tau\in K}\tau(c^*c)^{1/2}.
\]
The bounded tracial completion $C^K$ is the quotient of the
$\mathrm{C}^*$-algebra of bounded $\norm{\cdot}_{2,K}$-Cauchy sequences by the
ideal of bounded $\norm{\cdot}_{2,K}$-null sequences. Every trace in $K$
extends to this quotient. Its unit ball is the $\norm{\cdot}_{2,K}$
completion of the unit ball of $C$. We write $C^u=C^{\T(C)}$.
For a subalgebra $F\subseteq C$, $F^K$ means its bounded closure inside
$C^K$, using the restrictions of the designated traces.
For a displayed pair $F\subset C$, the notation $F^u$ always uses
the restrictions of $\T(C)$ and means the bounded closure in $C^u$.
\end{defn}

A nonempty compact trace family may be replaced by its closed convex
hull without changing the uniform $2$-seminorm or the completion.
This convention also permits compact designated subfamilies below.

Multiplication and involution are continuous in the uniform $2$-norm
on bounded sets, since
\begin{equation}\label{eq:two-mult}
 \norm{xay}_{2,K}\leq\norm{x}\norm{y}\norm{a}_{2,K}.
\end{equation}
The unit-ball description in Definition~\ref{def:completion} follows
by cutting bounded representatives down to the unit ball with
continuous functional calculus. All isomorphisms of completed pairs
in this paper preserve the designated trace sets. They are therefore
uniform $2$-norm isometries. An abstract $\mathrm{C}^*$-algebra isomorphism
without the trace marking is not the notion being classified.

We will also use bounded density in the following form. If a
$\mathrm{C}^*$-subalgebra $C_0\subset C$ is $2$-norm dense in $C$, its unit
ball is $2$-norm dense in the unit ball of $C$. Indeed, replace an
approximant $x\in C_0$ by its radial truncation
$r(x)=x\min\{1,|x|^{-1}\}$, defined continuously at zero. For a
contraction $a$, one has
$\norm{r(x)-a}_{2,\tau}\leq\norm{x-a}_{2,\tau}$ for each trace.
To see this in the finite von Neumann algebra of the trace, write
$x=v|x|$ and $t=(|x|-1)_+$. Then
$\operatorname{Re}\tau((x-r(x))^*(r(x)-a))
=\tau(t)-\operatorname{Re}\tau(tv^*a)\geq0$.
Expanding $\norm{x-a}_{2,\tau}^2$ proves the inequality. Taking the
supremum proves the assertion, also for completed ambient algebras.

\begin{defn}
A \emph{regular finite-dimensional map} into an AF algebra with its
standard diagonal sends the specified diagonal matrix units to
diagonal projections and every specified matrix unit to a partial
normalizer of the diagonal. Connecting maps of standard Bratteli
presentations are taken to be regular. A \emph{uniform Cartan pair}
$(\calD,\calM,K)$ consists of a tracially complete algebra, a MASA
$\calD$, a faithful trace-preserving conditional expectation onto it,
and normalizers whose linear span is dense in the uniform $2$-norm.
The last condition is a condition on the uniform topology.
\end{defn}

\subsection{Expectations and completed Cartan pairs}
\label{subsec:expectation-preliminaries}

\begin{defn}
A conditional expectation $\mathbb E:C\to F$, where $F$ is unital and
abelian, has \emph{Property (D)} if, for every finite $\calF\subset C$
and every $\eps>0$, there are contractions $x_1,\ldots,x_n\in F$ with
\[
 \norm{\mathbb E(a)-\frac1n\sum_{j=1}^n x_jax_j^*}<\eps
 \quad(a\in\calF).
\]
\end{defn}

\begin{lem}\label{lem:expectation}
If $\tau\circ\mathbb E=\tau$ for $\tau\in K$, then $\mathbb E$ is
$\norm{\cdot}_{2,K}$-contractive. Property (D) implies
$\tau\circ\mathbb E=\tau$ for every trace $\tau$ of $C$.
\end{lem}
\begin{proof}
Kadison's inequality gives
$\mathbb E(a)^*\mathbb E(a)\leq\mathbb E(a^*a)$; applying $\tau$ proves
the first assertion. For the second, apply Property (D) to $\{a,1\}$
with tolerances tending to zero. If
$q=n^{-1}\sum_jx_j^*x_j=n^{-1}\sum_jx_jx_j^*$, then $q\to1$ in norm,
and the trace of the average applied to $a$ is $\tau(aq)$.
\end{proof}

\begin{prop}\label{prop:uniform-cartan}
Suppose $F\subset C$ is a unital Cartan subalgebra and its expectation
has Property (D). Then $(F^K,C^K,K)$ is a uniform Cartan pair.
\end{prop}
\begin{proof}
Lemma~\ref{lem:expectation} extends $\mathbb E$ to a completely positive,
contractive, $F^K$-bimodular projection $C^K\to F^K$. This follows by
applying $\mathbb E$ to bounded Cauchy representatives; complete
positivity can equivalently be checked in the sequence-algebra
quotient. The extension preserves the designated traces and is
faithful, since $\mathbb E(a^*a)=0$ implies
$\norm{a}_{2,K}=0$.

For $a\in C^K$, approximate $a$ by a bounded element $b\in C$ and apply
Property (D) to $\{b,1\}$. Equation~\eqref{eq:two-mult} and
Lemma~\ref{lem:expectation} show that, with any prescribed accuracy,
\[
 \norm{\mathbb E(a)-\frac1n\sum_jx_jax_j^*}_{2,K}<\eps,
 \qquad \norm{1-q}_{2,K}<\eps.
\]
If $a$ commutes with $F$, the average equals $aq$, and therefore
$\norm{a-\mathbb E(a)}_{2,K}\leq(1+\norm{a})\eps$. Thus
$F'\cap C^K=F^K$, which proves the MASA assertion. The normalizers of
$F$ in $C$ still normalize $F^K$, by bounded $2$-norm continuity.
Cartan regularity in $C$ then implies uniform regularity in $C^K$.
\end{proof}

\begin{prop}\label{prop:free-averaging}
The canonical expectation for a free action has Property (D), with
the $x_j$ chosen to be unitaries. The standard AF expectation also
has Property (D).
\end{prop}
\begin{proof}
For the dynamical assertion it suffices to treat finitely many Fourier
polynomials. Let $L\subset\Gamma\setminus\{e\}$ contain their nonzero
Fourier indices. Freeness and compactness give finitely many real
functions $f_1,\ldots,f_r\in \mathrm{C}(X)$ and $c>0$ such that
\[
 \max_j|f_j(x)-f_j(g^{-1}x)|\geq c
 \quad(x\in X, g\in L).
\]
Average conjugation by
$v_t=\exp(i\sum_jt_jf_j)$ over $[-R,R]^r$, using normalized Lebesgue
measure. On the $g$th coefficient the multiplier is
\[
 \prod_{j=1}^r
 \frac{\sin(R(f_j(x)-f_j(g^{-1}x)))}
 {R(f_j(x)-f_j(g^{-1}x))},
\]
with value $1$ at zero denominators. Its modulus is at most $1/(Rc)$
for $g\in L$. The coefficient at $e$ is unchanged. Equal-weight
Riemann sums approximate this average in norm on the prescribed
finite set and consist of conjugations by unitaries in $D$.

For an AF algebra approximate the finite set in one finite stage.
Averaging diagonal phase unitaries at that stage gives exactly its
diagonal expectation. Regularity identifies this with the restriction
of the standard AF expectation.
\end{proof}

\section{Reduction to a Cantor system}\label{sec:cantor-reduction}

We use the zero-dimensional extension theorem of Kerr and Szab\'o
\cite[Theorem~5.5]{KS}: the small boundary property for an amenable
action is equivalent to the existence of a zero-dimensional compact
metrizable extension which is measure-isomorphic over singleton
fibres. In this statement the set of singleton fibres is conull for
every invariant measure, and the induced map on invariant-measure
simplexes is an affine homeomorphism.

\begin{prop}\label{prop:cantor-reduction}
Suppose that the action $\Gamma\curvearrowright X$ has the small boundary
property. There is a free minimal Cantor system $\Gamma\curvearrowright Y$ and an
equivariant surjection $\pi:Y\to X$ such that pullback gives a marked
isomorphism
\[
 (\mathrm{C}(X)^u,(\mathrm{C}(X)\rtimes_r\Gamma)^u)
 \cong(\mathrm{C}(Y)^u,(\mathrm{C}(Y)\rtimes_r\Gamma)^u).
\]
\end{prop}
\begin{proof}
Take the extension supplied by the cited theorem and restrict it to a
minimal closed invariant subset $Y$. Its image is $X$ by minimality.
Every singleton fibre lies in $Y$. Consequently all invariant lifts
of invariant measures on $X$ are supported on $Y$, and restriction
retains the measure-simplex homeomorphism and the singleton-fibre
property. Freeness is inherited. A minimal compact zero-dimensional
space for a free action of an infinite group has no isolated points,
and hence is a Cantor space.

The pullback crossed-product map is injective and preserves the
identified traces, so it extends isometrically to the bounded
completions. To prove diagonal density, let $V\subset Y$ be clopen.
Put
\[
 K_V=\pi(V)\cap\pi(Y\setminus V),\qquad
 U_V=X\setminus\pi(Y\setminus V).
\]
$K_V$ is closed and null for every invariant measure, and
$\partial U_V\subset K_V$. The indicator of $V$ agrees, outside
$\pi^{-1}(K_V)$, with the pullback of $1_{U_V}$. Continuous functions
on $X$ approximate $1_{U_V}$ uniformly in $2$-norm over the invariant
measures. Indeed, choose decreasing continuous contractions
$h_m\downarrow 1_{K_V}$ which are one on a neighborhood of $K_V$.
The functions $\mu\mapsto\mu(h_m)$ decrease to zero uniformly by
Dini's theorem. Each $(1-h_m)1_{U_V}$ is continuous and its squared
error is bounded by $h_m$. Locally constant
functions are norm dense in $\mathrm{C}(Y)$, proving diagonal density.
Fourier polynomials now prove density of the crossed products.
The radial truncation argument in Subsection~\ref{subsec:tracial-completions} makes this density hold
on unit balls, proving surjectivity of the completed map.
\end{proof}

\begin{lem}\label{lem:empirical}
Let $f_1,\ldots,f_r\in \mathrm{C}(Y)$ be real-valued and let $\eta>0$. There are a finite
$K\subset\Gamma$ and $\delta>0$ such that every nonempty finite
$(K,\delta)$-invariant set $S\subset\Gamma$ and every $y\in Y$ satisfy
\[
 \bigl(\abs{S}^{-1}\sum_{s\in S}f_j(sy)\bigr)_{j=1}^r
 \in\{(\mu(f_j))_{j=1}^r:\mu\in M_\Gamma(Y)\}+(-\eta,\eta)^r.
\]
Here $(K,\delta)$-invariance means
$|gS\mathbin\triangle S|<\delta|S|$ for $g\in K$.
\end{lem}
\begin{proof}
If the assertion fails, choose increasingly invariant finite sets
and offending points. A weak-star cluster point of their empirical
probability measures is invariant: the difference after translation
by any fixed $g$ is bounded by $|gS\triangle S|/|S|$ times the norm
of the test function. The prescribed finite coordinates converge to
the invariant-measure coordinates, a contradiction.
\end{proof}

\begin{lem}\label{lem:syndetic}
For every nonempty open $V\subset Y$ there are a finite set
$P\subset\Gamma$ and $c_V>0$ such that every sufficiently invariant
finite $S$ satisfies
$|\{s\in S:sy\in V\}|\geq c_V|S|$ for all $y\in Y$.
\end{lem}
\begin{proof}
Minimality gives $Y=\bigcup_{p\in P}p^{-1}V$. For each $s$ in the
$P$-interior of $S$, choose $p\in P$ with $psy\in V$. The map
$s\mapsto ps$ has fibres of size at most $|P|$. For sufficiently
invariant $S$ its interior has at least $|S|/2$ points, so take
$c_V=1/(2|P|)$.
\end{proof}

\section{A simple AF model with large multiplicities}\label{sec:AF-model}

For a finite $K\subset\Gamma$, put
$\operatorname{Int}_K(S)=\{s\in S:Ks\subset S\}$.
The clopen tower theorem \cite[Theorem~3.13]{KS} supplies, with any
prescribed invariance and coverage tolerances, a finite clopen castle
$\{(V_i,S_i)\}_i$: all levels $sV_i$ are disjoint, its shapes are
approximately invariant, and its complement has uniformly small
invariant measure. Its shapes can be required to exceed any prescribed
size, by strengthening the invariance requirement. For example, if
$|S|\leq N$ and $K$ has more than $N$ elements, sufficiently small
$K$-boundary is impossible. Splitting bases makes every level refine
any prescribed finite clopen partition. Reanchoring a tower permits
$e\in S_i$.

\begin{thm}\label{thm:AF-model}
Under the dynamical conventions of Section~1 and the small boundary
property, there are a simple unital infinite-dimensional AF algebra
$B$ with standard diagonal $E$ and a marked isomorphism
\[
 (D^u,A^u)\cong(E^u,B^u).
\]
Given any sequence of integers $R_n\geq2$, the model can be presented
by regular unital inclusions $B_n\to B_{n+1}$ whose every Bratteli
multiplicity is at least $R_n$.
\end{thm}

\begin{proof}
By Proposition~\ref{prop:cantor-reduction} we work on a Cantor space
$Y$. We give the finite-stage construction and then verify both
inverse maps. All $2$-norms on $A$ use its full trace simplex; a
$2$-norm on a finite-dimensional algebra uses all its tracial states.

\medskip\noindent\emph{Orbit compression and physical lifts.}
For a castle $(V_i,S_i)$ choose $z_i\in V_i$, set $p_i=|S_i|$, and put
\[
 B_n=\bigoplus_i\mathrm{M}_{p_i},\qquad E_n=\bigoplus_i\C^{p_i}.
\]
The orbit representation on $\ell^2(\Gamma)$ is
$\rho_z(f)\delta_s=f(sz)\delta_s$ and
$\rho_z(u_g)\delta_s=\delta_{gs}$. The compressions
\[
 \phi_n(a)=\bigoplus_i P_{S_i}\rho_{z_i}(a)P_{S_i}
\]
are unital completely positive contractions and send $D$ into $E_n$.
There is an exact, generally nonunital homomorphism
\[
 L_n:B_n\longrightarrow A,\qquad
 L_n(e^i_{s,t})=u_s1_{V_i}u_t^*\quad(s,t\in S_i).
\]
Castle disjointness proves the matrix-unit relations. Its matrix
units normalize $D$, and
\begin{equation}\label{eq:lift-left-inverse}
 \phi_nL_n=\id_{B_n},\qquad
 q_n:=L_n(1)=1_{\bigcup_iS_iV_i}.
\end{equation}
Write $\delta_n=\sup_\mu\mu(1-q_n)$.

If $a=\sum_gf_gu_g$, all $f_g$ have oscillation at most $\omega$ on
each level, and $|S_i\setminus g^{-1}S_i|/|S_i|\leq\eta$, then
\begin{equation}\label{eq:physical-approx}
 \norm{L_n\phi_n(a)-a}_{2,\T(A)}
 \leq\sum_g\bigl(\omega+\norm{f_g}\sqrt{\eta+\delta_n}\bigr).
\end{equation}
Indeed, the lifted compression of $u_g$ is $u_g$ times the indicator
of source levels whose $g$-translate is in the same tower. The omitted
source set has measure at most $\eta+\delta_n$ for every invariant
measure; the coefficient error on the retained levels is at most
$\omega$. Also, in any block with compression projection $P$,
\[
 \phi_n(ab)-\phi_n(a)\phi_n(b)
 =P\rho(a)(1-P)\rho(b)P.
\]
For $b=\sum_hb_hu_h$, its normalized Hilbert--Schmidt norm is at most
\[
 \norm{a}\sum_h\norm{b_h}\sqrt{|S\setminus h^{-1}S|/|S|}.
\]
These estimates are uniform over all block traces and control any
fixed finite Fourier family and its products.

\medskip\noindent\emph{Complete copies of the old towers.}
Fix castle $n$, and consider a later tower $(W,T)$. Split $W$ by all
the clopen membership tests $s^{-1}tV_i^{(k)}$, for $s\in T$,
$t\in S_i^{(k)}$, and $k\leq n$. Thus all earlier castle membership
patterns are constant on its levels. The complete old tiles
\[
 S_i h\subset T,\qquad hy\in V_i\quad(y\in W)
\]
and their types are now independent of $y$. They are disjoint:
freeness turns equality of orbit coordinates into equality of group
coordinates, and the old castle levels are disjoint.

Let $C_n=Y\setminus\bigcup_iS_iV_i$ and
$K_n=\bigcup_iS_iS_i^{-1}$. Every coordinate in
$\operatorname{Int}_{K_n}(T)$ whose orbit point is outside $C_n$
belongs to a complete old tile. Lemma~\ref{lem:empirical}, applied to
$1_{C_n}$, bounds the fraction of unmatched coordinates by
$2\delta_n+\zeta_n$, after sufficiently strong new invariance.
Here $\zeta_n>0$ is arbitrarily small.

Every old type occurs arbitrarily often in every new tower.
For clarity, choose finite $F_i$ with $Y=\bigcup_{f\in F_i}fV_i$.
For each $q\in\operatorname{Int}_{F_i^{-1}\cup S_iF_i^{-1}}(T)$,
some $f\in F_i$ gives $hy\in V_i$ for $h=f^{-1}q$, and
$S_i h\subset T$. Each $h$ is counted at most $|F_i|$ times.
Consequently the number $m_i$ of complete type-$i$ copies is at least
\[
 |\operatorname{Int}_{F_i^{-1}\cup S_iF_i^{-1}}(T)|/|F_i|.
\]
Arbitrarily large and sufficiently invariant shapes make this exceed
any prescribed constant, uniformly over the new towers.

\medskip\noindent\emph{A unital arithmetic repair.}
Put $d=\gcd(p_1,\ldots,p_r)$. Choose a conductor $C\geq1$ so that
every integer at least $C$ is a nonnegative integer combination of
$p_1/d,\ldots,p_r/d$, and set $Q=dC$. Trim fewer than $d$ coordinates
from each new shape so that its size $q$ is divisible by $d$.
This loses at most $d/\min q$ of coverage, and changes every fixed
boundary proportion by a quantity tending to zero as $\min q\to\infty$.
Apply the preceding packing estimates to these trimmed shapes.

Require $m_1\geq R_n+\lceil Q/p_1\rceil$ and $m_i\geq R_n$ for
$i>1$. The residual rank $r=q-\sum_i m_ip_i$ is divisible by $d$.
If $r=0$ or $r\geq Q$, it is representable by the old block sizes.
If $0<r<Q$, delete the least number of complete type-$1$ copies
needed to bring the residual rank to at least $Q$. This deletes at
most $\lceil Q/p_1\rceil$ copies, so every retained multiplicity
is still at least $R_n$. The deleted rank is at most $Q+p_1$.
Partition the residual coordinates into blocks of the old sizes
and place arbitrary coordinate matrix-unit copies on them.

Doing this in each new block defines a regular unital injective
homomorphism $\psi_n:B_n\to B_{n+1}$. All its multiplicities are at
least $R_n$. It agrees with the natural orbit identification on the
retained complete tiles. Call the other coordinates exceptional.
Their fraction in every new block is at most
\begin{equation}\label{eq:exceptional-fraction}
 f_n=2\delta_n+\zeta_n+(Q+p_1)/\min_j|T_j|.
\end{equation}
The old remainder term is essential. The other terms, including the
trimming loss absorbed into $\zeta_n$, can be made arbitrarily small
once the old castle has been fixed.

\medskip\noindent\emph{Compatibility estimates.}
On retained copies the old and new coefficient evaluations differ
by the oscillation on an old level. Away from old shape boundaries,
their group permutation edges agree. Thus, for a fixed Fourier
polynomial $a=\sum_gf_gu_g$, an absolute constant, say $4$, gives
\begin{equation}\label{eq:forward-bound}
 \norm{\phi_{n+1}(a)-\psi_n\phi_n(a)}_{2,\T(B_{n+1})}
 \leq4\sum_g\bigl(\omega_n(f_g)+\norm{f_g}\sqrt{\eta_n+f_n}\bigr).
\end{equation}
Here $\eta_n$ bounds the relevant old shape boundary proportions.
For the term $f_gu_g$, call a source coordinate good if it belongs
to a retained tile $S_i h$ and its $g$-translate belongs to that
same tile. On good coordinates both matrices use the same translation
and their coefficient values differ by at most $\omega_n(f_g)$.
The remaining columns form a fraction at most $f_n+\eta_n$ of the
block. On their span the difference has operator norm at most
$2\norm{f_g}$, so its normalized Hilbert--Schmidt norm is at most
$2\norm{f_g}\sqrt{f_n+\eta_n}$. Summing over $g$ gives
\eqref{eq:forward-bound}.

There is a sharper estimate on the whole newly created lift ball:
\begin{equation}\label{eq:new-lift-bound}
 \norm{\phi_{n+1}L_n(b)-\psi_n(b)}_{2,\T(B_{n+1})}
 \leq2\norm{b}\sqrt{f_n}.
\end{equation}
Indeed, on a retained tile $S_i h$ over a new base $W$, the old
matrix unit $L_n(e^i_{s,t})$ sends the coordinate $th$ to $sh$ in
the orbit representation sampled on $W$. It vanishes on the other
columns of that tile and on retained tiles of other types, by
disjointness of the old levels. This is exactly the copy used by
$\psi_n$. The same calculation for adjoints shows that the retained
coordinate projection reduces both maps. On its complement both are
contractions. In the physical algebra the same matrix unit carries
$thW$ onto $shW$ by the translation $st^{-1}$, as does its new copy.
Thus the indicator $r_n$ of all retained tiles reduces both
$L_n(B_n)$ and $L_{n+1}\psi_n(B_n)$, which agree on $r_n$. Since invariant measures give equal measure to every level
of a tower, $\mu(1-r_n)\leq\delta_{n+1}+f_n$. Hence
\begin{equation}\label{eq:reverse-bound}
 \norm{L_{n+1}\psi_n(b)-L_n(b)}_{2,\T(A)}
 \leq2\norm{b}\sqrt{\delta_{n+1}+f_n}.
\end{equation}

\medskip\noindent\emph{Scheduling.}
Fix a summable sequence $\eps_n>0$, for example $2^{-n-20}$.
Before choosing castle $n$, choose a finite Fourier test family
$\calF_n$. Fix a norm-dense sequence of Fourier polynomials in
the unit ball of $A$. At step $n$ include the first $n$ polynomials,
their adjoints and pairwise products, all previous tests, and a finite
$\eps_n/8$-net of $L_i((B_i)_1)$ for every $i<n$. These lift balls
are norm-compact sets of Fourier polynomials. The castle at stage $n$ makes the
physical and multiplicative errors on these tests at most $\eps_n$.
Its old boundaries and level oscillations are also made small enough
for the outgoing estimate \eqref{eq:forward-bound}. For example, if
$M_n$ is the largest sum of coefficient norms in the scheduled
Fourier tests, require
$2\delta_n<(\eps_n/(100(1+M_n)))^2$, and budget their old boundary
and oscillation errors below $\eps_n/4$.

When choosing castle $n+1$, the old dimensions, conductor, membership
tests, and all constants in those estimates are already fixed.
The empirical lemma, the complete-tile count, and arbitrarily large
shapes make the remaining terms of \eqref{eq:exceptional-fraction}
small enough that \eqref{eq:new-lift-bound} and
\eqref{eq:reverse-bound} are at most $\eps_n\norm{b}$ and
\eqref{eq:forward-bound} is at most $\eps_n/2$ on $\calF_n$ and on
the selected lift nets. Take their meshes at most $\eps_n/8$.
The two maps are norm contractions, so the net errors add at most
twice the mesh, giving a bound $\eps_n$ on every previous lift
unit ball. Each new lift ball has its own direct estimate
\eqref{eq:new-lift-bound}. All conditions at a stage are finite and
are satisfied by strengthening the invariance, refining the clopen
partition, and decreasing the coverage error. This is an induction,
with no condition on a previously chosen castle introduced later.

\medskip\noindent\emph{Limits and inverse identities.}
Let $B=\varinjlim(B_n,\psi_n)$, with inclusions $\iota_n$, and
$E=\overline{\bigcup_n\iota_n(E_n)}$. The maps are unital and all
incidences are positive; consequently $B$ is simple and unital.
Since $R_n\geq2$, it is infinite-dimensional. $E$ is its standard
AF Cartan subalgebra and $\T(B)$ is nonempty and compact.
Summable forward errors and norm density define a unital homomorphism
\[
 \Phi:A\to B^u,\qquad \Phi(a)=\lim_n\iota_n\phi_n(a).
\]
The multiplicative estimates prove multiplicativity, and
$\Phi(D)\subset E^u$. Summable reverse errors define a homomorphism
\[
 \beta:B\to A^u,\qquad
 \beta(\iota_n(b))=\lim_{m\to\infty}L_m\psi_{n,m}(b),
\]
where $\psi_{n,m}$ denotes the connecting map. It is unital because
$L_m(1)=q_m\to1$ in uniform $2$-norm, and $\beta(E)\subset D^u$.

No finite-stage trace contraction is asserted. Instead, every
$\tau\in\T(B)$ pulls back along the unital homomorphism $\Phi$ to
a trace of $A$, and every $\rho\in\T(A)$ pulls back along $\beta$
to a trace of $B$. Therefore
\[
 \norm{\Phi(a)}_{2,\T(B)}\leq\norm{a}_{2,\T(A)},\qquad
 \norm{\beta(b)}_{2,\T(A)}\leq\norm{b}_{2,\T(B)}.
\]
Both maps extend to the bounded completions.

Put $t_n=\sum_{k\geq n}\eps_k\to0$. The scheduled lift estimates
and \eqref{eq:lift-left-inverse} prove, uniformly for $b\in(B_n)_1$,
\begin{equation}\label{eq:two-tail-estimates}
 \norm{\Phi L_n(b)-\iota_n(b)}_2\leq t_n,
 \qquad \norm{\beta\iota_n(b)-L_n(b)}_2\leq t_n.
\end{equation}
For the first estimate, fix $b\in(B_n)_1$. For every $m\geq n$,
\[
 \norm{\phi_{m+1}(L_n b)-\psi_m\phi_m(L_n b)}_{2,\T(B_{m+1})}
 \leq\eps_m.
\]
At $m=n$ this is \eqref{eq:new-lift-bound}, since $\phi_nL_n(b)=b$.
For $m>n$ it follows from the net for the fixed lift ball
$L_n((B_n)_1)$ chosen before castle $m$. Each trace of $B$ restricts
to a trace of $B_{m+1}$, so these bounds telescope after the stage
inclusions. The second estimate telescopes
\eqref{eq:reverse-bound} along the connecting maps.
For an element of a fixed $B$ stage, view it in arbitrarily late
stages and apply \eqref{eq:two-tail-estimates} and contractivity:
$\norm{\Phi\beta(b)-b}_2\leq2t_n\norm{b}\to0$.
For a fixed scheduled contraction $a$, physical approximation gives
\[
 \norm{\beta\Phi(a)-a}_2
 \leq t_n+t_n+\norm{L_n\phi_n(a)-a}_2
 \leq2t_n+\eps_n\longrightarrow0.
\]
Density proves the two inverse identities on the full completions.
They preserve the completed diagonals in both directions.

Finally, the trace pullbacks are inverse affine maps between the
\emph{full} trace simplexes. They are continuous because, for a
fixed element of either original algebra, its image is a uniform
$2$-norm limit of fixed elements of the other algebra; trace
evaluation is then a uniform limit of continuous affine functions.
Thus the maps are affine homeomorphisms, and $\Phi$ is exactly
uniform-$2$-isometric with this trace marking. It also intertwines
the diagonal expectations. Indeed, the transported expectation and
the canonical expectation both preserve every designated trace.
For any $x$, their difference $d$ at $x$ is diagonal. Bimodularity
and trace preservation give $\tau(d^*d)=0$ for every designated
trace, since both terms in $\tau(d^*d)$ equal $\tau(d^*x)$.
Faithfulness of the designated family gives $d=0$.
\end{proof}

\begin{rem}
The finite compressions are generally not $2$-contractive relative
to $\T(A)$: a continuous function equal to one at all sampled points
may have very small invariant-measure $2$-norm. The limiting unital
homomorphisms, followed by their inverse identities, supply the
required trace control. The theorem constructs an isomorphism with
the AF \emph{completion}; it asserts no exact diagonal-preserving
embedding $A\to B$.
\end{rem}

\section{Relative classification of the completed AF pairs}\label{sec:AF-classification}

All AF algebras in this section are separable, simple, unital, and
infinite-dimensional, with standard diagonals and regular unital
finite-dimensional presentations. We use the standard dimension-group
facts
\begin{equation}\label{eq:dimension-group}
 \overline{\rho_C(\mathrm{K}_0(C))}=\Aff(\T(C)),\qquad
 \mathrm{K}_0(C)^+=\{0\}\cup\{g:\tau(g)>0\text{ for every }\tau\in\T(C)\}.
\end{equation}
Here $\rho_C(g)(\tau)=\tau(g)$, and the AF trace simplex is the
normalized state space $\T(C)=S(\mathrm{K}_0(C),[1_C])$. See, for example,
\cite[Theorem~2.6 and Proposition~2.8]{Dadarlat}.
Positive classes, their comparisons, and equalities are witnessed by
rank vectors at sufficiently late AF stages. In particular a positive
class below $[1]$ is represented by a diagonal projection, and equal
classes of diagonal projections can be matched by a finite-stage
partial permutation. These observations will retain the diagonals.

\begin{lem}\label{lem:regular-existence}
Let $F=\bigoplus_{i=1}^r\mathrm{M}_{n_i}$ and let
$f_i\in\Aff(\T(C))$ be strictly positive, with
$\sum_i n_if_i\leq1$. For every $\eta>0$ there is an injective regular
finite-stage homomorphism $\alpha:F\to C$, possibly nonunital, such
that
\[
 \left|\tau(\alpha(x))-\sum_i f_i(\tau)\Tr_{n_i}(x_i)\right|
 \leq\eta\norm{x}\quad(\tau\in\T(C),\ x\in F).
\]
\end{lem}
\begin{proof}
Put $L=\sum_i n_i$ and choose $\gamma>0$ with
$\gamma<\frac12\min_{i,\tau}f_i(\tau)$ and
$\frac54L\gamma<\eta$. By \eqref{eq:dimension-group}, choose classes
$g_i$ with
$\norm{\rho_C(g_i)-(f_i-\gamma)}_\infty<\gamma/4$.
They are positive and
$[1]-\sum_i n_ig_i$ is positive, since its value at every trace is
at least $\frac34L\gamma$.
At a sufficiently late stage partition the diagonal unit into
$n_i$ projections of class $g_i$ for each $i$, and a residual
projection. The rank-vector description of the equality
$[1]=\sum_i n_ig_i+([1]-\sum_i n_ig_i)$ justifies this partition.
For each $i$, coordinate bijections between its equal-rank projections
give matrix units. The resulting homomorphism is injective and
regular. Off-diagonal matrix units have zero trace, and
$|\Tr(x_i)|\leq n_i\norm{x}$, giving the asserted estimate.
\end{proof}

\begin{lem}\label{lem:regular-uniqueness}
Let $\alpha,\beta:F\to C$ be regular finite-stage homomorphisms,
where $F=\bigoplus_i\mathrm{M}_{n_i}$ and $L=\sum_i n_i$. Suppose $\alpha$
is injective and put
$c=\min_{i,\tau}\tau(\alpha(e^i_{11}))>0$.
If $0<\delta<c/4$ and
\[
 |\tau(\alpha(e^i_{11}))-\tau(\beta(e^i_{11}))|<\delta
 \quad(i,\tau),
\]
there is a finite-stage unitary $U\in C$ normalizing its standard
diagonal such that
\[
 \sup_{x\in F_1}\norm{U\alpha(x)U^*-\beta(x)}_{2,\T(C)}
 <4\sqrt{L\delta}.
\]
\end{lem}
\begin{proof}
Write $p_i=\alpha(e^i_{11})$ and $q_i=\beta(e^i_{11})$.
Choose $r_i\in \mathrm{K}_0(C)$ within $\delta/4$, on all traces, of
$\rho_C([p_i])-2\delta$. Then
\[
 \tfrac74\delta<\tau(p_i)-\tau(r_i)<\tfrac94\delta,\qquad
 \tfrac34\delta<\tau(q_i)-\tau(r_i)<\tfrac{13}4\delta.
\]
Also $\tau(r_i)>c-\frac94\delta>0$. Strict order yields
$0<r_i<[p_i],[q_i]$. At a common later stage choose diagonal
subprojections $p'_i\leq p_i$, $q'_i\leq q_i$ of class $r_i$.
Transport these along the respective regular matrix units, obtaining
\[
 P=\sum_{i,k}\alpha(e^i_{k1})p'_i\alpha(e^i_{1k}),\qquad
 Q=\sum_{i,k}\beta(e^i_{k1})q'_i\beta(e^i_{1k}).
\]
These diagonal projections reduce the corresponding representations.
Choose partial permutations $w_i$ from $p'_i$ to $q'_i$. Then
\[
 w=\sum_{i,k}\beta(e^i_{k1})w_i\alpha(e^i_{1k})
\]
has initial projection $P$, final projection $Q$, and intertwines the
cut representations. The summands have orthogonal diagonal initial
and final supports, so $w$ is a partial normalizer. Since
$[1-P]=[1-Q]$, a later-stage partial permutation $w_0$ matches
their complements. Thus $U=w+w_0$ is a unitary normalizer.

If $E=\alpha(1_F)$ and $H=\beta(1_F)$, the remaining discrepancy is
\[
 U\alpha(x)(E-P)U^*-\beta(x)(H-Q).
\]
For $\norm{x}\leq1$, its $2$-norm is at most
$\sqrt{\frac94L\delta}+\sqrt{\frac{13}4L\delta}<4\sqrt{L\delta}$.
This also explains why the finite maps need not be unital.
\end{proof}

\begin{thm}\label{thm:AF-pairs}
Let $(E_C,C)$ and $(E_B,B)$ be standard pairs of simple unital
infinite-dimen\-sional AF algebras. Every affine homeomorphism
$h:\T(B)\to\T(C)$ is implemented by an isomorphism
\[
 \Phi:(E_C^u,C^u)\longrightarrow(E_B^u,B^u),\qquad
 \tau(\Phi(x))=h(\tau)(x).
\]
It intertwines the diagonal expectations.
\end{thm}
\begin{proof}
We give the summable intertwining, including the $2$-norm transport
estimates. Choose exhausting subsequences $F_n\subset C$,
$G_n\subset B$ of the given unital stages and regular homomorphisms
\[
 \alpha_n:F_n\to G_n,\qquad \beta_n:G_n\to F_{n+1}.
\]
Write $i_n,j_n$ for the stage inclusions. Fix $\eps_n=2^{-n}$.
We arrange
\begin{align}
 \sup_{x\in(F_n)_1}\norm{\beta_n\alpha_n(x)-i_n(x)}_{2,\T(C)}
 &<\eps_n,\label{eq:AF-first-diagram}\\
 \sup_{y\in(G_n)_1}\norm{\alpha_{n+1}\beta_n(y)-j_n(y)}_{2,\T(B)}
 &<\eps_n,\label{eq:AF-second-diagram}
\end{align}
and trace errors $a_n,b_n\leq\eps_n^2$ on whole stage unit balls:
\begin{align*}
 |\tau(\alpha_n(x))-h(\tau)(x)|&\leq a_n\norm{x},\\
 |\sigma(\beta_n(y))-h^{-1}(\sigma)(y)|&\leq b_n\norm{y}.
\end{align*}

Here is the tolerance induction. For a fixed stage inclusion
$F\subset C$ or $G\subset B$, set
\[
 \Delta_F(\eps)=\min\left\{
 \tfrac14\min_{i,\tau}\tau(e^i_{11}),\,
 \frac{\eps^2}{16\sum_i n_i}\right\}.
\]
Choose $a_n$ smaller than $\eps_n^2$ and
$\frac14\Delta_{F_n}(\eps_n)$, and, for $n>1$, smaller than
$\frac14\Delta_{G_{n-1}}(\eps_{n-1})$. The previous step will have
chosen $b_{n-1}<\frac14\Delta_{G_{n-1}}(\eps_{n-1})$.
Lemma~\ref{lem:regular-existence}, with
$f_i(\tau)=h(\tau)(e^i_{11})$, gives $\alpha_n$ with this error.
For $n>1$ the minimal-projection trace error of
$\alpha_n\beta_{n-1}$ against inclusion is at most $a_n+b_{n-1}$.
Lemma~\ref{lem:regular-uniqueness} conjugates $\alpha_n$ by a
diagonal-normalizing unitary so that
\eqref{eq:AF-second-diagram} holds for $n-1$.
Choose $G_n$ late enough to contain the corrected image, the previous
stage, the finite-stage correction, and the $n$th original $B$ stage.

Now choose $b_n$ smaller than $\eps_n^2$ and both
$\frac14\Delta_{F_n}(\eps_n)$ and
$\frac14\Delta_{G_n}(\eps_n)$. Apply existence in $C$ with the
transported trace profiles of $G_n$. The trace error of
$\beta_n\alpha_n$ against inclusion is at most $a_n+b_n$.
Relative uniqueness permits a diagonal-normalizing correction of
$\beta_n$ giving \eqref{eq:AF-first-diagram}. Choose $F_{n+1}$
to contain all these data and the next original $C$ stage.
All corrections preserve the trace estimates and regularity.
The minima are positive by simplicity, faithfulness of traces, and
compactness; thus this induction is possible.

Applied to $z^*z$, the trace estimates imply
\begin{equation}\label{eq:transport-two}
 \norm{\alpha_n(z)}_{2,\T(B)}^2
 \leq\norm{z}_{2,\T(C)}^2+a_n\norm{z}^2,
 \qquad
 \norm{\beta_n(z)}_{2,\T(C)}^2
 \leq\norm{z}_{2,\T(B)}^2+b_n\norm{z}^2.
\end{equation}
Insert $\alpha_{n+1}\beta_n\alpha_n(x)$ between
$\alpha_{n+1}i_n(x)$ and $j_n\alpha_n(x)$.
Equations~\eqref{eq:AF-first-diagram}--\eqref{eq:transport-two} give
uniformly for $x\in(F_n)_1$
\[
 \norm{\alpha_{n+1}i_n(x)-j_n\alpha_n(x)}_2
 \leq\sqrt{\eps_n^2+4a_{n+1}}+\eps_n
 \leq2\eps_n+2\sqrt{a_{n+1}}=:c_n.
\]
Similarly, uniformly for $y\in(G_n)_1$,
\[
 \norm{\beta_{n+1}j_n(y)-i_{n+1}\beta_n(y)}_2
 \leq\eps_n+\eps_{n+1}+2\sqrt{b_{n+1}}=:d_n.
\]
Both sequences are summable. The factors $4$ account for the
operator-norm bound $2$ on the differences being transported.

Hence on the algebraic stage unions the limits
$\Phi(x)=\lim_n\alpha_n(x)$ and $\Psi(y)=\lim_n\beta_n(y)$ exist
in the bounded completions. They are contractive homomorphisms by
bounded continuity of multiplication. Their limiting trace identities
are exact, since $a_n,b_n\to0$. Applying them to $x^*x$ gives
exact uniform-$2$-isometries. These maps therefore extend to the
original algebras and their bounded completions.

For $C_n=\sum_{r\geq n}c_r$ and $D_n=\sum_{r\geq n}d_r$, the
limits satisfy
$\norm{\Phi(x)-\alpha_n(x)}_2\leq C_n$ on $(F_n)_1$ and
$\norm{\Psi(y)-\beta_n(y)}_2\leq D_n$ on $(G_n)_1$.
For a contraction in a fixed $F_k$, at every $n\geq k$,
\[
 \norm{\Psi\Phi(x)-x}_2\leq C_n+D_n+\eps_n\longrightarrow0.
\]
For a contraction in a fixed $G_k$, the symmetric estimate is
$\norm{\Phi\Psi(y)-y}_2\leq D_n+C_{n+1}+\eps_n\to0$.
Thus the maps are inverse on the full completions, and in particular
unital. Every finite map preserves the diagonals, so the inverse
limits carry the completed diagonals onto each other. A regular
finite map intertwines the diagonal expectations: its off-diagonal
matrix units have disjoint diagonal initial and final supports, and
thus zero diagonal expectation. Passing to the limits proves the
last assertion.
\end{proof}

\begin{proof}[Proof of Theorem~\ref{thm:main}]
Apply Theorem~\ref{thm:AF-model} to both actions. Compose the induced
affine trace identifications with the specified $H$, and apply
Theorem~\ref{thm:AF-pairs} to their AF models, using the inverse
orientation for $h$. Composition gives the stated marked pair
isomorphism. Conversely the trace pullback of a marked isomorphism
is affine, bijective, and continuous: evaluations on the original
algebras are uniform limits of continuous affine evaluations on
bounded approximants. Its inverse has the same property.
\end{proof}

\begin{proof}[Proof of Remark~\ref{rem:intro-Z-stable}]
Under URP, \cite[Theorem~1.2]{ENrc} gives
$\operatorname{rc}(\mathrm{C}(X)\rtimes_r\Gamma)
\geq\tfrac12\operatorname{mdim}(X,\Gamma)$.
$\mathcal Z$-stability gives radius of comparison zero, so the action has
zero mean dimension. URP then gives the small boundary property by
\cite[Theorem~5.1]{NiuZ}.
\end{proof}

\begin{rem}
Simplicity and infinite dimension are substantive in the AF pair
theorem. $\mathrm{M}_2$ and $\mathrm{M}_3$ have the same singleton trace simplex but
different completed pairs. Even the cross maps between infinite
simple AF models need not be unital at finite stages, as a unital
map $\mathrm{M}_2\to \mathrm{M}_{3^\infty}$ does not exist. Their missing corners
vanish in the tracial limit.
\end{rem}

\begin{lem}\label{lem:large-telescoping}
Any Bratteli presentation by unital inclusions of a simple
infinite-dimensional AF algebra admits a telescoping with prescribed lower bounds on
every multiplicity.
\end{lem}
\begin{proof}
The minimum size of its late matrix blocks tends to infinity.
For any $N$, choose $N$ orthogonal nonzero projections in a finite
stage. Each is full, so for each such projection $p$ there is
$M_p$ with $[1]\leq M_p[p]$. These inequalities become rank-vector
inequalities in all sufficiently late stages. Thus every late block
contains every one of these $N$ orthogonal nonzero images.

Now fix a coordinate minimal projection $p_i$ in each summand of
the source stage. Fullness gives $[1]\leq M_i[p_i]$. In a late
target block of size $s_j$ and source multiplicity $m_{ji}$ this
says $s_j\leq M_im_{ji}$. The minimum $s_j$ tends to infinity,
so all $m_{ji}$ exceed any prescribed bound at sufficiently late
stages. Choose those stages successively.
\end{proof}

\section{Borel realization and orbit equivalence}\label{sec:uniform-orbits}

We describe the dynamical meaning of an isomorphism of marked uniform tracial completions.
Throughout this section actions are left actions, with
\(\alpha_g(f)(x)=f(g^{-1}x)\) and \(u_gfu_g^*=\alpha_g(f)\).

Let \(X\) be compact metrizable and let \(\Delta\) be a nonempty compact set of Borel probability measures on \(X\).
Put
\[
 \norm{f}_{2,\Delta}=\sup_{\mu\in\Delta}\left(\int_X\abs{f}^2\,d\mu\right)^{1/2},
 \qquad \mathcal D_\Delta=\overline{\mathrm{C}(X)}^{\Delta}.
\]
The completion is taken on norm-bounded sets, after quotienting by the seminorm kernel if necessary.
A Borel set is \emph{\(\Delta\)-null} if it has measure zero for every \(\mu\in\Delta\).
In this section, ``universally conull'' always refers to the designated measure family.
Write \(\mathcal B_\Delta(X)\) for the algebra of bounded Borel functions modulo functions vanishing outside a \(\Delta\)-null set.
Its norm is \(\sup_{\mu\in\Delta}\operatorname*{ess\,sup}_{\mu}\abs{b}\).

\begin{lem}\label{lem:orbit-borel-realization}
There is a canonical faithful unital homomorphism
\[
 \mathcal D_\Delta\longrightarrow\mathcal B_\Delta(X).
\]
Its image consists of the bounded Borel classes uniformly approximable in \(\norm{\cdot}_{2,\Delta}\) by a norm-bounded sequence of continuous functions.
For every \(b\in\mathcal D_\Delta\), the functions
\(\mu\mapsto\int b\,d\mu\) and \(\mu\mapsto\int\abs{b}^2\,d\mu\) are continuous on \(\Delta\).
Moreover, \(\mathcal D_\Delta\) embeds faithfully into \(L^\infty(X,\mu_0)\), where
\(\mu_0=\sum_{n\geq1}2^{-n}\mu_n\) and \((\mu_n)\) is any weak*-dense sequence in \(\Delta\).
\end{lem}

\begin{proof}
Let \((f_n)\) be a uniform \(2\)-Cauchy sequence with \(\norm{f_n}\leq M\).
Choose a subsequence satisfying
\(\sum_j\norm{f_{n_{j+1}}-f_{n_j}}_{2,\Delta}<\infty\).
For every \(\mu\in\Delta\),
\[
 \int_X\sum_j\abs{f_{n_{j+1}}-f_{n_j}}\,d\mu
 \leq\sum_j\norm{f_{n_{j+1}}-f_{n_j}}_{2,\Delta}<\infty.
\]
Thus the subsequence converges pointwise outside one Borel \(\Delta\)-null set.
Its limit \(b\), extended by zero on that set, is Borel and bounded by \(M\).
Minkowski's inequality and Fatou's lemma give
\[
 \norm{b-f_{n_j}}_{2,\Delta}
 \leq\sum_{k\geq j}\norm{f_{n_{k+1}}-f_{n_k}}_{2,\Delta}\longrightarrow0.
\]
Equivalent sequences give limits agreeing almost everywhere for every measure in \(\Delta\).
Addition and involution pass to these limits; multiplication does also, since on bounded sets
\[
 \norm{f_ng_n-fg}_{2,\Delta}
 \leq\norm{f_n}\norm{g_n-g}_{2,\Delta}+\norm{g}\norm{f_n-f}_{2,\Delta}.
\]
This proves the homomorphism assertion, faithfulness, and the image description.
The integrals of \(f_n\) and \(\abs{f_n}^2\) converge uniformly on \(\Delta\) to the stated integrals of \(b\), proving continuity.
For the last assertion, passing to \(L^\infty(X,\mu_0)\) is well-defined because every \(\Delta\)-null set is \(\mu_0\)-null.
If \(\int\abs{b}^2\,d\mu_0=0\), then this integral vanishes for every \(\mu_n\), and continuity makes it vanish on all of \(\Delta\).
Hence \(b=0\) in \(\mathcal D_\Delta\).
\end{proof}

\begin{rem}
The measure \(\mu_0\) need not dominate every measure in \(\Delta\), and a \(\mu_0\)-null Borel set need not be \(\Delta\)-null.
Faithfulness in the lemma concerns the completed diagonal.
We use \(\mathcal B_\Delta(X)\) for statements involving all measures simultaneously.
\end{rem}

\begin{thm}\label{thm:orbit-spatial-diagonal}
Let \(X,Y\) be compact metrizable, and let \(\Delta_X,\Delta_Y\) be nonempty compact sets of probability measures.
Suppose \(H:\Delta_X\to\Delta_Y\) is a homeomorphism and
\(\Phi:\mathcal D_{\Delta_X}\to\mathcal D_{\Delta_Y}\) is a unital isomorphism such that
\[
 \int_Y\Phi(b)\,d(H\mu)=\int_Xb\,d\mu
 \quad(\mu\in\Delta_X,\ b\in\mathcal D_{\Delta_X}).
\]
Then there are Borel universally conull sets \(X^0\subseteq X\), \(Y^0\subseteq Y\), and a Borel isomorphism \(\theta:X^0\to Y^0\), such that
\(\theta_*\mu=H\mu\) for every \(\mu\in\Delta_X\) and
\[
 \Phi([b])=[b\circ\theta^{-1}]\quad\text{in }\mathcal B_{\Delta_Y}(Y).
\]
Functions on the conull domains may be extended by zero to their complements.
\end{thm}

\begin{proof}
The trace correspondence makes \(\Phi\) an isometry for the uniform \(2\)-norms.
Choose a countable norm-dense unital *-subalgebra \(\mathcal A\subseteq\mathrm{C}(X)\), over \(\mathbb Q+i\mathbb Q\), containing a separating family.
Choose a Borel representative \(b_a\) of \(\Phi(a)\), bounded by \(\norm{a}\), for each \(a\in\mathcal A\).
All rational linear, multiplicative, involutive, and unital relations between these representatives hold outside a universally null set.
There are countably many relations, so they hold simultaneously on one Borel universally conull set \(S_Y\).
At each \(y\in S_Y\), evaluation \(a\mapsto b_a(y)\) extends to a character of \(\mathrm{C}(X)\).
It is therefore evaluation at a unique point \(\sigma(y)\in X\).
A countable continuous separating family shows that \(\sigma\) is Borel.
Extend it arbitrarily to the omitted set.
The trace identity and norm density imply
\(\sigma_*(H\mu)=\mu\), by equality of integrals of continuous functions.

If bounded continuous \(a_n\) approximate a Borel representative \(b\) of an arbitrary completed element, then
\[
 \sup_{\nu\in\Delta_Y}\int_Y\abs{a_n\circ\sigma-b\circ\sigma}^2\,d\nu
 =\sup_{\mu\in\Delta_X}\int_X\abs{a_n-b}^2\,d\mu\longrightarrow0.
\]
Consequently \(\Phi([b])=[b\circ\sigma]\) for all completed \(b\).
Apply the same argument to \(\Phi^{-1}\), obtaining a Borel map \(\rho:X\to Y\), a good set \(S_X\), and \(\rho_*\mu=H\mu\).
Applying the two compositions to countably many separating functions gives
\(\sigma\rho(x)=x\) universally almost everywhere and \(\rho\sigma(y)=y\) universally almost everywhere.
Preimages of the respective universally null sets are universally null, by the pushforward identities.
Set
\[
 \begin{split}
 X^0&=S_X\cap\rho^{-1}(S_Y)\cap\{x:\sigma\rho(x)=x\},\\
 Y^0&=S_Y\cap\sigma^{-1}(S_X)\cap\{y:\rho\sigma(y)=y\}.
 \end{split}
\]
These are Borel universally conull sets, and \(\theta=\rho|_{X^0}\) is a Borel bijection onto \(Y^0\), with inverse \(\sigma|_{Y^0}\).
\end{proof}

\begin{rem}
The pullback identity is an equality of Borel classes.
For independently chosen representatives its pointwise exceptional set may depend on the completed element.
Exact pointwise equality on the fixed domain is obtained by defining image representatives by pullback, without intersecting uncountably many exceptional sets.
\end{rem}

Now let a countable group \(\Gamma\) act freely on \(X\), let
\(K\subseteq M_\Gamma(X)\) be nonempty and compact, and write
\(\mathcal M=\overline{\mathrm{C}(X)\rtimes_r\Gamma}^{K}\), \(\mathcal D=\overline{\mathrm{C}(X)}^K\).
The canonical expectation \(\mathbb E\) extends to
\(\mathbb E:\mathcal M\to\mathcal D\), since it is contractive for every tracial \(2\)-norm.
For finite \(F\subseteq\Gamma\), define
\(a_h(z)=\mathbb E(zu_h^*)\) and \(P_F(z)=\sum_{h\in F}a_h(z)u_h\).
For each \(\mu\in K\), \(P_F\) is an orthogonal projection in \(L^2(\tau_\mu)\), so
\[
 \norm{z-P_Fz}_{2,\tau_\mu}^2
 =\norm{z}_{2,\tau_\mu}^2-\sum_{h\in F}\int\abs{a_h(z)}^2\,d\mu.
\]
Moreover \(P_Fz\to z\) in the uniform \(2\)-norm as finite sets exhaust \(\Gamma\): for a finite Fourier polynomial \(b\) supported in \(F\),
\(\norm{z-P_Fz}_{2,K}\leq2\norm{z-b}_{2,K}\).

\begin{lem}\label{lem:orbit-uniform-normalizer}
If \(v\in\mathcal M\) is unitary and \(v\mathcal Dv^*=\mathcal D\), its diagonal automorphism is induced by a Borel bijection \(T\) between universally conull subsets of \(X\), preserving every \(\mu\in K\).
Its graph lies in the orbit relation.
The source sets \(Q_h=\{x:T(x)=hx\}\) have indicators \(q_h\in\mathcal D\), and
\[
 \inf_{F\subseteq\Gamma\ \mathrm{finite}}\sup_{\mu\in K}
 \mu\left(X\setminus\bigcup_{h\in F}Q_h\right)=0.
\]
\end{lem}

\begin{proof}
Conjugation \(\beta(d)=vdv^*\) preserves every trace, so Theorem~\ref{thm:orbit-spatial-diagonal} supplies \(T\) with \(\beta(d)=[d\circ T^{-1}]\).
Writing \(a_h=\mathbb E(vu_h^*)\), the equation \(vf=\beta(f)v\) gives
\(a_h\alpha_h(f)=\beta(f)a_h\) for \(f\in\mathrm{C}(X)\).
On the nonzero support of \(a_h\), a countable separating family therefore gives \(T^{-1}y=h^{-1}y\).
Freeness makes the supports for distinct \(h\) disjoint.
These conclusions hold simultaneously outside one universally null set.
Since \(\norm{a_h}\leq1\), the pointwise sum \(\sum_h\abs{a_h}^2\) is at most one.
Parseval says its integral is one for every \(\mu\in K\), so the sum equals one universally almost everywhere.
Disjointness now makes each \(e_h=\abs{a_h}^2\) a projection of \(\mathcal D\).
The Fourier coefficients of \(v^*\), indexed by \(h^{-1}\), are \(\alpha_{h^{-1}}(a_h^*)\).
Applying the same argument to \(v^*\) shows that \(q_h=\alpha_{h^{-1}}(e_h)\) are mutually orthogonal projections partitioning one.
Their sets are exactly the source sets of \(T(x)=hx\).
Finally invariance of the measures and the Fourier projection estimate give
\[
 \sup_{\mu\in K}\sum_{h\notin F}\int q_h\,d\mu
 =\norm{v-P_Fv}_{2,K}^2\longrightarrow0.
\]
Remove the $\Gamma$-saturation of the universally null complements
of the good source and range sets, and of the exceptional graph
identities for $T$ and $T^{-1}$. The remaining Borel set is invariant
and universally conull. Since both maps are orbit maps there, it is
preserved by $T$ and $T^{-1}$.
\end{proof}

\begin{defn}\label{defn:uniform-measure-oe}
Let countable amenable groups \(\Gamma_i\) act freely on compact metrizable \(X_i\), and let \(K_i=M_{\Gamma_i}(X_i)\neq\varnothing\).
A \emph{uniform measure orbit equivalence} implementing an affine homeomorphism \(H:K_1\to K_2\) is a Borel bijection \(\theta:X_1^0\to X_2^0\) between invariant universally conull Borel sets such that:
\begin{enumerate}
\item \(\theta(\Gamma_1x)=\Gamma_2\theta(x)\) and \(\theta_*\mu=H\mu\);
\item \(f\circ\theta^{-1}\in\mathcal D_2\) for \(f\in\mathrm{C}(X_1)\), and \(k\circ\theta\in\mathcal D_1\) for \(k\in\mathrm{C}(X_2)\);
\item every partition \(\{x:\theta(gx)=h\theta(x)\}\) has indicator in \(\mathcal D_1\), and every analogous partition for \(\theta^{-1}\) has indicator in \(\mathcal D_2\).
\end{enumerate}
Here \(\mathcal D_i=\overline{\mathrm{C}(X_i)}^{K_i}\). Membership means uniform $L^2(K_i)$ approximability by a bounded sequence of continuous functions, so all conditions refer to the spaces and measure families.
\end{defn}

For a fixed \(g\), the cocycle \(c(g,x)\) defined by \(\theta(gx)=c(g,x)\theta(x)\) has uniformly small tails over \(K_1\).
Indeed the finite sums
\[
 \mu\longmapsto\sum_{h\in F}\mu\{c(g,\cdot)=h\}
\]
are continuous, increase to one, and converge uniformly by Dini's theorem.
This imposes no bound on cocycle word lengths; the terminology is distinguished from bounded orbit equivalence.

\subsection{The uniform measure full group}

For the full-group statements, let $\Gamma$ be countable and act freely
on $X$, and let $K\subseteq M_\Gamma(X)$ be nonempty and compact.
Put $\mathcal D=\mathrm{C}(X)^K$, $\mathcal M=(\mathrm{C}(X)\rtimes_r\Gamma)^K$, and
$R_\Gamma=\{(x,gx):x\in X,\ g\in\Gamma\}$.

\begin{defn}\label{defn:uniform-measure-full-group}
The \emph{uniform measure full group} $[R_\Gamma]^u_K$ consists of
Borel bijections $T:X^0\to X^0$ of invariant universally conull Borel
sets, with $T(x)\in\Gamma x$, such that, for every $g\in\Gamma$,
$\mathbf1_{\{x:T(x)=gx\}}$ is uniformly approximable in $L^2(K)$ by
continuous functions taking values in $[0,1]$. Equivalently,
\[
 q_{T,g}:=\mathbf1_{\{x:T(x)=gx\}}\in\mathcal D\quad(g\in\Gamma).
\]
Maps agreeing almost everywhere for every measure in $K$ are
identified. Extending by the identity off $X^0$ gives an equivalent
formulation using Borel bijections of $X$. Define
\[
 d_K(T,S)=\sup_{\mu\in K}\mu\{x:T(x)\ne S(x)\}.
\]
All sets and functions are understood modulo simultaneous null sets.
\end{defn}

\begin{prop}\label{prop:uniform-full-group-complete}
The set $[R_\Gamma]^u_K$ is a group, and $d_K$ is a
complete bi-invariant metric. Every element preserves every measure in
$K$ and induces a uniform $2$-isometric automorphism of $\mathcal D$ by
$d\mapsto[d\circ T^{-1}]$.
\end{prop}

\begin{proof}
Write $Q_g=\{T=g\}$. The sets $Q_g$ partition the domain, and the sets
$gQ_g$ partition the range. For an arbitrary Borel set $B$, invariance
of $\mu$ under the original action gives
\[
 \mu(T^{-1}B)=\sum_g\mu(Q_g\cap g^{-1}B)
 =\sum_g\mu(gQ_g\cap B)=\mu(B).
\]
The finite sums $\sum_{g\in F}\mu(Q_g)$ are continuous on $K$, increase
to one, and converge uniformly by Dini's theorem. Consequently both
\[
 \sum_{g\in F}q_{T,g}\alpha_{g^{-1}}(f),
 \qquad
 \sum_{g\in F}\alpha_g(q_{T,g}f)
\]
converge uniformly in $2$-norm, with a common operator-norm bound
$\norm f$. Their Borel limits are respectively $f\circ T$ and
$f\circ T^{-1}$. Hence these pullbacks belong to $\mathcal D$ for
$f\in\mathrm{C}(X)$, and extend to inverse uniform $2$-isometries of
$\mathcal D$.

The inverse partitions have indicators
$q_{T^{-1},g}=\alpha_{g^{-1}}(q_{T,g^{-1}})$. For a product $TS$, its
partition for the arrow $k$ is the disjoint union of
\[
 Q_{S,g}\cap g^{-1}Q_{T,h},\qquad hg=k.
\]
Their indicators belong to $\mathcal D$. Keeping only finitely many
$g$ and $h$ omits a set of measure at most the corresponding tail of
$S$ plus that of $T$, because $S$ preserves every measure. These tails
are uniformly small. Thus every product partition indicator is a
bounded uniform $2$-limit in $\mathcal D$, proving closure under
composition.

The triangle inequality follows from inclusion of disagreement sets.
Left multiplication preserves these sets since the multiplying map is
injective, and right multiplication preserves their measures. Therefore
$d_K$ is bi-invariant; in particular
$d_K(T^{-1},S^{-1})=d_K(T,S)$.

To prove completeness, choose a subsequence $(T_n)$ with
\[
 \sum_n d_K(T_n,T_{n+1})<\infty.
\]
On one invariant
universally conull Borel set all the maps and their inverses are
defined. The Borel sets
\[
 \limsup_n\{T_n\ne T_{n+1}\},\qquad
 \limsup_n\{T_n^{-1}\ne T_{n+1}^{-1}\}
\]
are universally null, by the summable estimates and the
Borel--Cantelli argument applied to each measure. Remove their
$\Gamma$-saturations. The remaining set $X^0$ is invariant and
universally conull, and every $T_n$ and $T_n^{-1}$ preserves $X^0$.
For each $x\in X^0$ the two sequences of point values eventually
stabilize. Their limits $T(x)$ and $S(x)$ are Borel orbit maps and
are inverse: once $T_n(x)=T(x)$, the equality
$T_n^{-1}(T(x))=x$ and stabilization at $T(x)$ imply $S(T(x))=x$.
The other composition is identical.

For this subsequence,
\[
 d_K(T_n,T)\leq\sum_{j\geq n}d_K(T_j,T_{j+1})\longrightarrow0.
\]
For each $g$,
\[
 \norm{q_{T_n,g}-\mathbf1_{\{T=g\}}}_{2,K}^2
 \leq d_K(T_n,T).
\]
Closure of the bounded uniform $2$-completion shows that the limit
indicator belongs to $\mathcal D$. Thus $T$ lies in the full group.
The original Cauchy sequence has the same limit.
\end{proof}

\begin{prop}\label{prop:uniform-full-group-normalizers}
Let $\mathcal N_{\mathcal M}(\mathcal D)$ be the unitary normalizer of
$\mathcal D$ in $\mathcal M$. There is a split exact sequence
\[
 1\longrightarrow\mathcal U(\mathcal D)
 \longrightarrow\mathcal N_{\mathcal M}(\mathcal D)
 \longrightarrow[R_\Gamma]^u_K
 \longrightarrow1.
\]
The canonical splitting is
\[
 T\longmapsto v_T:=\sum_{g\in\Gamma}u_gq_{T,g},
\]
where the sums converge on norm-bounded sets in the uniform $2$-norm.
It satisfies
\[
 \norm{v_T-v_S}_{2,K}^2=2d_K(T,S).
\]
If normalizers $v,v'$ induce $T,S$, respectively, then
\[
 \inf_{w\in\mathcal U(\mathcal D)}\norm{v-wv'}_{2,K}^2
 =2d_K(T,S).
\]
Thus the full group is precisely the normalizer modulo diagonal phase
unitaries, with its quotient topology. In particular
$\mathcal N_{\mathcal M}(\mathcal D)$ is the semidirect product
$\mathcal U(\mathcal D)\rtimes[R_\Gamma]^u_K$.
\end{prop}

\begin{proof}
The finite sums defining $v_T$ are contractions: their source
projections and range projections are pairwise orthogonal. The
uniform tail estimate in the preceding proof makes them uniformly
$2$-Cauchy. Their limit is a unitary and implements
$d\mapsto[d\circ T^{-1}]$.
For products, the identity
\[
 (u_hq_{T,h})(u_gq_{S,g})
 =u_{hg}\alpha_{g^{-1}}(q_{T,h})q_{S,g}
\]
and the product partitions show that $v_Tv_S=v_{TS}$, first on
finite pieces and then by bounded uniform $2$-continuity.
Lemma~\ref{lem:orbit-uniform-normalizer} gives surjectivity. If $v$
induces $T$, then $v_T^*v$ commutes with $\mathcal D$ and hence belongs
to $\mathcal U(\mathcal D)$ by Propositions~\ref{prop:uniform-cartan}
and~\ref{prop:free-averaging}.
This also identifies the kernel and the unique diagonal phase.

The identity arrow of $T^{-1}S$ is exactly the agreement set of $T$
and $S$. Freeness of the action therefore gives
\[
 \mathbb E(v_T^*v_S)=e_{T,S}:=\mathbf1_{\{T=S\}}\in\mathcal D,
 \qquad
 \norm{v_T-v_S}_{2,\tau_\mu}^2=2(1-\mu\{T=S\}).
\]
Taking suprema proves the first metric identity.
Write $v=v_Ta$ and $v'=v_Sb$ with $a,b\in\mathcal U(\mathcal D)$.
For right multiplication by a diagonal unitary $c$,
\[
 \mathbb E(v^*v'c)=a^*e_{T,S}bc.
\]
Its real trace is at most $\mu\{T=S\}$ for every $\mu$. Thus every
choice of $c$ gives a squared uniform distance at least $2d_K(T,S)$.
Taking $c=b^*a$ gives equality, by unitary invariance of the
$2$-norm. Since $v'$ normalizes $\mathcal D$, left multiplication
$wv'$ ranges over the same set as right multiplication $v'c$.
This proves the quotient metric formula. Finally, with range phases
$v=wv_T$, multiplication is
$(w,T)(z,S)=(w[z\circ T^{-1}],TS)$, which is the stated
semidirect product. The quotient map is continuous by the metric
formula, and its canonical section is continuous by the first metric
identity. This makes the splitting topological. The unit ball of
$\mathcal M$ is separable for the uniform $2$-norm, since the crossed
product is norm separable. The canonical section therefore makes the
full group separable as well. Together with completeness this shows
that $[R_\Gamma]^u_K$ is a Polish group.
\end{proof}

\begin{rem}
For $K=\{\mu\}$ the indicator requirement is automatic, by the density
of continuous functions in $L^2(\mu)$. The definition recovers the
usual measured full group and its uniform metric; see
\cite{CLM} for this classical setting. For a family of measures,
uniform continuous approximability of the individual indicators is
substantive. Uniformly small displacement tails alone do not imply it.
No amenability, minimality, or small boundary property is needed for
the preceding full-group results.
\end{rem}

\begin{cor}\label{cor:uniform-oe-full-group}
If $\theta$ is a uniform measure orbit equivalence implementing
$H:K_1\to K_2$, then
\[
 T\longmapsto\theta T\theta^{-1}
\]
is an isometric isomorphism
$[R_{\Gamma_1}]^u_{K_1}\cong[R_{\Gamma_2}]^u_{K_2}$.
\end{cor}

\begin{proof}
For fixed $h\in\Gamma_2$, pullback by $\theta$ of the conjugated
displacement indicator is
\[
 \mathbf1_{Q_h(\theta T\theta^{-1})}\circ\theta
 =\sum_{g\in\Gamma_1}\mathbf1_{Q_g(T)}\,
 \mathbf1_{\{x:\theta(gx)=h\theta(x)\}}.
\]
Finite subsums belong to $\mathcal D_1$, and their error is
supported on the uniformly small displacement tail of $T$.
Thus this indicator belongs to $\mathcal D_1$. The diagonal
pullback isomorphism carries it into $\mathcal D_2$.
The inverse construction uses $\theta^{-1}$.
Finally $\theta_*\mu=H\mu$ identifies the disagreement sets and
their measures, so it preserves the metric.
\end{proof}

The topological full group $[[\Gamma\curvearrowright X]]$ consists of
homeomorphisms whose orbit cocycle is continuous into the discrete
group $\Gamma$, and hence has finite range. For earlier work on
uniform topologies and topological-full-group density in Cantor
minimal systems, see Bezuglyi--Kwiatkowski~\cite{BK}.

\begin{prop}[Density in the Cantor case]\label{prop:uniform-full-group-density}
Let a countably infinite amenable group $\Gamma$ act freely on a compact
zero-dimensional metrizable space $X$, and put $K=M_\Gamma(X)$.
Then the topological full group $[[\Gamma\curvearrowright X]]$ is dense
in $[R_\Gamma]^u_K$. Indeed, every $T\in[R_\Gamma]^u_K$ is a
$d_K$-limit of finite-order elements of the topological full group.
Consequently $[R_\Gamma]^u_K$ is the metric completion of the
countable topological full group, modulo simultaneous-null equality.
\end{prop}

\begin{proof}
Write $Q_h=\{x:T(x)=hx\}$. These are a countable partition of a
universally conull set. Since their indicators belong to $\mathcal D$,
the functions $\mu\mapsto\mu(Q_h)$ are continuous on $K$. Dini's theorem
therefore supplies a finite set $F\subseteq\Gamma$ with
\[
 \sup_{\mu\in K}\mu\left(X\setminus\bigcup_{h\in F}Q_h\right)
 <\varepsilon/4.
\]
Every projection of $\mathcal D$ is uniformly approximable in measure
by clopen sets: approximate its indicator in uniform $2$-norm by a
continuous function in $[0,1]$, approximate that function uniformly by
a locally constant function, and threshold at $1/2$.
Thus choose clopen $C_h$, $h\in F$, so that, with
$Z=\bigcup_{h\in F}(C_h\mathbin\triangle Q_h)$,
\[
 \sup_{\mu\in K}\mu(Z)<\varepsilon/8.
\]

We first turn these candidate arrows into an injective clopen partial
map. Order $F$ and put
\[
 P_h=C_h\setminus\bigcup_{k<h}C_k,\qquad
 U_h=P_h\setminus h^{-1}\left(\bigcup_{k<h}kP_k\right).
\]
The $U_h$ are clopen and pairwise disjoint, and their translates $hU_h$
are pairwise disjoint. The partial map $x\mapsto hx$ on $U_h$ agrees
with $T$ on all but uniformly less than
$\varepsilon/4+2\sup_\mu\mu(Z)<\varepsilon/2$ of $X$.
Here is the measure estimate. Outside $Z$, membership in $P_h$ is
exactly membership in $Q_h$. If a correct arrow from $P_h\cap Q_h$ is
removed because its range meets an earlier $kP_k$, that earlier arrow
has an incorrect source: otherwise injectivity of $T$ and freeness
would force the two arrows to be identical. The union of the ranges
of incorrect arrows has measure at most $\mu(Z)$, since the $P_k$
are disjoint and the group elements preserve $\mu$. Pulling it back by
the measure-preserving map $T$ bounds the additional loss by
$\mu(Z)$.

By almost finiteness in measure in the zero-dimensional case
\cite[Theorem~3.13]{KS}, choose a clopen castle
$\{(V_i,S_i)\}_i$ with uniformly small remainder and shapes sufficiently
$F$-invariant that
\[
 \sup_{\mu\in K}\mu\left(X\setminus\bigcup_iS_iV_i\right)
 +\sum_{h\in F}\sum_i
       \abs{\{s\in S_i:hs\notin S_i\}}\mu(V_i)
 <\varepsilon/2.
\]
Retain an arrow with source in $U_h\cap sV_i$ only when $hs\in S_i$.
These retained arrows form an injective partial map within each
individual finite tower. Refine each base $V_i$ by the finitely many
clopen membership patterns $s^{-1}U_h$. On a resulting base atom $W$,
the retained arrows give a fixed partial injection from $S_i$ to
$S_i$. Complete it to a permutation of $S_i$. On the corresponding
levels, define the homeomorphism by
\[
 sy\longmapsto \sigma_{i,W}(s)y\qquad(y\in W),
\]
and use the identity outside the castle footprint.
The resulting homeomorphism $S$ has a locally constant finite-valued
$\Gamma$-cocycle, so $S\in[[\Gamma\curvearrowright X]]$. It has finite
order because it acts by finitely many finite permutations. The two
preceding loss estimates show $d_K(T,S)<\varepsilon$.
Finally, a compact metrizable zero-dimensional space has only
countably many clopen sets, so its topological full group is countable.
\end{proof}

\begin{rem}
This density assertion need not hold for the topological full group
on the original nonzero-dimensional space. For an irrational
rotation $R$ of the circle, every continuous $\mathbb Z$-valued cocycle
is constant. Thus $[[\mathbb Z\curvearrowright\mathbb T]]=
\{R^n:n\in\mathbb Z\}$, and different powers have uniform-measure
distance one. If an interval $I$ is disjoint from $R(I)$, the Borel
involution which is $R$ on $I$, $R^{-1}$ on $R(I)$, and the identity
elsewhere belongs to the uniform-measure full group but lies a positive
distance from all these powers. For general SBP systems, the Cantor
extension gives a useful zero-dimensional model; it does not identify
its topological full group with the topological full group on the
original space.
\end{rem}

\begin{cor}\label{cor:canonical-normalizer-density}
In the Cantor setting of Proposition~\ref{prop:uniform-full-group-density},
every canonical orbit unitary $v_T$ is a uniform $2$-limit of
finite-order unitary normalizers belonging to the original crossed
product. Their orbit maps lie in $[[\Gamma\curvearrowright X]]$.
\end{cor}
\begin{proof}
For a topological-full-group element $S$, the finite sum
$v_S=\sum_g u_gq_{S,g}$ belongs to the crossed product, since the
$q_{S,g}$ are clopen indicators. The canonical section is a
homomorphism, so finite order of $S$ implies finite order of $v_S$.
Now use the density proposition and
$\norm{v_T-v_S}_{2,K}^2=2d_K(T,S)$.
\end{proof}

\subsection{Finite permutation groups from the AF model}

\begin{lem}\label{lem:AF-normalizer-density}
Let $(E,B)$ be the standard diagonal pair of a simple unital AF algebra,
with coordinate Bratteli inclusions $B_n\subseteq B_{n+1}$.
Put $\mathcal D=E^u$ and $\mathcal M=B^u$.
Let $P_n$ be the finite group of coordinate permutation unitaries
in $B_n$. For every
$v\in\mathcal N_{\mathcal U(\mathcal M)}(\mathcal D)$ and every
$\eta>0$, there are $m$, $p\in P_m$, and
$\lambda\in\mathcal U(\mathcal D)$ such that
\[
 \norm{v-\lambda p}_{2,\T(B)}<\eta.
\]
The groups $P_n$ are increasing, and their union is dense in the
completed normalizer quotient modulo diagonal unitaries.
\end{lem}

\begin{proof}
Write $B_n=\bigoplus_{a=1}^s\mathrm{M}_{k_a}$, with coordinate matrix units
$e^a_{ij}$ and central projections $z_a=\sum_i e^a_{ii}$.
Put $\mathcal D_a=e^a_{11}\mathcal D$ and
\[
 C_n=\mathrm{C}^*(B_n,\mathcal D)
 \cong\bigoplus_{a=1}^s\mathrm{M}_{k_a}(\mathcal D_a).
\]
The diagonal in this matrix description is $\mathcal D$.
If $\mathbb E:\mathcal M\to\mathcal D$ is the extended diagonal expectation,
define
\[
 \mathsf E_n(x)=
 \sum_{a=1}^s\sum_{i,j=1}^{k_a}
 e^a_{i1}\,\mathbb E(e^a_{1i}xe^a_{j1})\,e^a_{1j}.
\]
To see that this is a conditional expectation onto $C_n$, first
pinch $x$ into the corners $z_a\mathcal Mz_a$, identify each corner
with $\mathrm{M}_{k_a}(e^a_{11}\mathcal Me^a_{11})$, and apply the corner
expectation $\mathbb E$ entrywise. This is unital and completely positive,
fixes $C_n$, and is $C_n$-bimodular.
For every designated trace $\tau$, cyclicity shows
\[
 \tau(\mathsf E_n(x))
 =\sum_{a,i}\tau(\mathbb E(e^a_{1i}xe^a_{i1}))
 =\sum_{a,i}\tau(xe^a_{ii})=\tau(x).
\]
Thus $\mathsf E_n$ is the $L^2(\tau)$-orthogonal projection onto
the closure of $C_n$, and is uniformly $2$-contractive.
Since it fixes $B_n$ and $\bigcup_n B_n$ is uniformly $2$-dense,
\[
 \norm{v-\mathsf E_n(v)}_{2,\T(B)}\longrightarrow0.
\]

Put $b=\mathsf E_n(v)$ and $\beta=\operatorname{Ad}(v)|_{\mathcal D}$.
The $\mathcal D$-bimodule property gives $bd=\beta(d)b$.
Consequently $r=v^*b$ commutes with $\mathcal D$.
The completed diagonal is maximal abelian, so $r\in\mathcal D$.
Idempotence now gives
\[
 b=\mathsf E_n(b)=\mathsf E_n(vr)=\mathsf E_n(v)r=vr^2,
\]
and hence $r^2=r$. An idempotent in a commutative $\mathrm{C}^*$-algebra
is a projection. Thus
\[
 b=vr,\qquad
 b^*b=r,\qquad bb^*=\beta(r),\qquad
 \norm{v-b}_{2,\tau}^2=\tau(1-r).
\]
In particular $b$ is an exact partial normalizer of $\mathcal D$
inside $C_n$.

We next complete $b$ to a monomial unitary in $C_n$.
In a block $\mathrm{M}_{k_a}(\mathcal D_a)$ write its entries as $b_{ij}$.
The diagonal conditions for $b e_{jj}b^*$ and $b^*e_{ii}b$
give, respectively,
\[
 b_{ij}b_{\ell j}^*=0\quad(i\ne\ell),\qquad
 b_{ij}^*b_{i\ell}=0\quad(j\ne\ell).
\]
The positives $p_{ij}=|b_{ij}|^2$ therefore have orthogonal supports
in every row and column. Their column sums are the diagonal entries
of the projection $r$, so each $p_{ij}$ is a projection.
The finite Boolean algebra generated by these projections partitions
the corner unit $e^a_{11}$ into projections on which the nonzero
entries have one fixed partial-permutation pattern.
Extend every such pattern to a permutation by matching the unused
columns and rows in increasing order. Give each added entry phase
one, and retain the existing entries of $b$.
This produces a monomial unitary $w\in C_n$ with $wr=b$.
Put $c=w-b=w(1-r)$.
Then $c\in C_n$, $c^*c=1-r$, and
\[
 \tau((v-b)^*c)
 =\tau\bigl(\mathsf E_n(v-b)^*c\bigr)=0.
\]
It follows that
\begin{equation}\label{eq:AF-partial-normalizer-rounding}
 \norm{v-w}_{2,\tau}^2=2\tau(1-r).
\end{equation}
In each block its row phase is explicitly
\[
 \lambda_i=\sum_j b_{ij}+\left(1-\sum_j p_{ij}\right)
 \in\mathcal U(\mathcal D_a).
\]
Orthogonality of the row supports proves that these are unitaries.
Removing them writes $w=\lambda\sigma$, where
$\lambda\in\mathcal U(\mathcal D)$ and $\sigma\in C_n$ is a
permutation matrix over a finite projection partition in each
$\mathcal D_a$.

Here are explicit bounds for approximating $\sigma$ by a finite-stage
coordinate permutation. Let
$(q_{a,\ell})_{\ell=1}^{L_a}$ be the tail projection partition
selecting its permutations in block $a$.
At a common stage $m\geq n$, choose positive diagonal contractions
$f_{a,\ell}\in e^a_{11}E_m e^a_{11}$ approximating
$q_{a,\ell}$ uniformly in $2$-norm.
For $\ell<L_a$ put
\[
 p_{a,\ell}=\mathbf1_{[1/2,1]}(f_{a,\ell}).
\]
These are finite-stage diagonal projections and satisfy, pointwise,
\[
 |p_{a,\ell}-q_{a,\ell}|^2
 \leq4|f_{a,\ell}-q_{a,\ell}|^2.
\]
Disjointize the $p_{a,\ell}$ in order and put the remaining
projection into the last part. The resulting partition
$(\widehat q_{a,\ell})$ agrees with $(q_{a,\ell})$ outside
\[
 t_a=\bigvee_{\ell<L_a}|p_{a,\ell}-q_{a,\ell}|,
 \qquad
 \tau(t_a)\leq
 4\sum_{\ell<L_a}\norm{f_{a,\ell}-q_{a,\ell}}_{2,\tau}^2.
\]
Replace the selectors of $\sigma$ by these finite-stage selectors,
obtaining $p\in P_m$. This is a coordinate permutation because
the coordinate Bratteli maps copy each old matrix unit along identical
tail coordinates. The two unitaries agree off the central projection
$\widetilde t_a=\sum_i e^a_{i1}t_ae^a_{1i}$, and
$\tau(\widetilde t_a)=k_a\tau(t_a)$. Hence
\begin{equation}\label{eq:AF-selector-rounding}
 \norm{\sigma-p}_{2,\T(B)}^2
 \leq16\sum_{a=1}^s k_a\sum_{\ell<L_a}
 \norm{f_{a,\ell}-q_{a,\ell}}_{2,\T(B)}^2.
\end{equation}
Choose $n$ so late that
$\norm{v-w}_{2,\T(B)}<\eta/2$, using
\eqref{eq:AF-partial-normalizer-rounding}.
There are only finitely many selectors, so choose their approximants
to make \eqref{eq:AF-selector-rounding} less than $\eta^2/4$.
Then
\[
 \norm{v-\lambda p}_{2,\T(B)}
 \leq\norm{v-w}_{2,\T(B)}+\norm{\sigma-p}_{2,\T(B)}<\eta.
\]
Finally coordinate Bratteli inclusions send permutation matrices to
permutation matrices, so $P_n\subseteq P_{n+1}$.
Every finitely generated subgroup of $\bigcup_nP_n$ is contained in
one finite $P_n$, proving local finiteness as well as density in
the normalizer quotient.
\end{proof}

\begin{cor}\label{cor:uniform-full-group-locally-finite}
For every dynamical system covered by Theorem~\ref{thm:main}, its
uniform measure full group has a dense locally finite subgroup.
Consequently its finite-order elements are dense.
\end{cor}

\begin{proof}
Transport the finite-stage AF permutation groups through the
completed-pair isomorphism in Theorem~\ref{thm:AF-model}, and then
pass from normalizers to their point transformations using
Proposition~\ref{prop:uniform-full-group-normalizers}.
Their images are increasing finite groups, so their union is locally
finite. If $v_T$ is a graph unitary, transport it to the AF completion
and apply Lemma~\ref{lem:AF-normalizer-density}.
Diagonal phases do not affect the point transformation, and
\[
 2d_K(T,S)\leq\norm{v_T-\lambda v_S}_{2,K}^2.
\]
The transported finite-group images are therefore dense for $d_K$.
Every element in their union has finite order.
\end{proof}

\begin{rem}
The dense locally finite subgroup is supplied by the AF model.
Its elements need not belong to the topological full group of the
original system: the original space may be connected, and the
model isomorphism controls completed diagonal projections rather
than clopen subsets of that space.
\end{rem}

\subsection{Orbit equivalence and completed pairs}

\begin{thm}\label{thm:uniform-measure-oe}
With the notation of Definition~\ref{defn:uniform-measure-oe}, put
\[
 \mathcal M_i=\overline{\mathrm{C}(X_i)\rtimes_r\Gamma_i}^{K_i}.
\]
There is an isomorphism \(\Phi:\mathcal M_1\to\mathcal M_2\) satisfying
\(\Phi(\mathcal D_1)=\mathcal D_2\) and
\(\tau_{H\mu}(\Phi(z))=\tau_\mu(z)\)
if and only if there is a uniform measure orbit equivalence implementing \(H\).
\end{thm}

\begin{proof}
Suppose first that \(\Phi\) exists.
Theorem~\ref{thm:orbit-spatial-diagonal} supplies a Borel bijection \(\theta:Y_1\to Y_2\) between universally conull sets, with the required trace transport and diagonal pullbacks.
For \(g\in\Gamma_1\), \(\Phi(u_g)\) normalizes \(\mathcal D_2\).
Its point map agrees with \(\theta g\theta^{-1}\) universally almost everywhere, by applying the diagonal conjugation identity to countably many separating coordinates.
Lemma~\ref{lem:orbit-uniform-normalizer} gives its orbit arrows and source partition projections in \(\mathcal D_2\); pulling these back by \(\Phi^{-1}\) gives the required projections in \(\mathcal D_1\).
The argument for \(\Phi^{-1}\) gives the reverse arrows and partitions.
All these assertions are countable and may therefore be made simultaneous on a Borel universally conull good set \(A_0\subseteq Y_1\).
Include in \(A_0\) the conditions \(gx\in Y_1\) and \(h\theta(x)\in Y_2\) for every \(g,h\), together with all forward and reverse orbit identities.
They are universally conull by invariance and trace transport.
Set \(X_1^0=\bigcap_g g^{-1}A_0\) and \(X_2^0=\theta(X_1^0)\).
The first is invariant and universally conull; the second is Borel and universally conull.
If \(x\in X_1^0\), the reverse identity gives \(\theta^{-1}(h\theta(x))=kx\) for some \(k\in\Gamma_1\).
Thus \(X_2^0\) is invariant as well, and \(\theta\) maps entire orbits onto entire orbits there.
Phase factors of the normalizers do not affect this conclusion.

Conversely, suppose that \(\theta\) is a uniform measure orbit equivalence.
The observable conditions and trace transport extend pullback to a uniform \(2\)-isometric isomorphism
\(\Phi_D:\mathcal D_1\to\mathcal D_2\).
Put \(T_g=\theta g\theta^{-1}\) and \(q_{g,h}=\mathbf1_{\{y:T_g(y)=hy\}}\in\mathcal D_2\).
The source projections \(q_{g,h}\) partition one; their range projections \(\alpha_h(q_{g,h})\) also partition one, since \(T_g\) is bijective.
Thus $T_g$ belongs to the uniform measure full group on $X_2$.
Proposition~\ref{prop:uniform-full-group-normalizers} gives its canonical
unitary $v_g=v_{T_g}$, and its homomorphic section gives
$v_gv_k=v_{gk}$ without a phase obstruction. The diagonal action of
that unitary is precisely the covariance relation
\(v_g\Phi_D(f)v_g^*=\Phi_D(\alpha_g(f))\).
The universal property and amenability give a homomorphism
\(\Psi:\mathrm{C}(X_1)\rtimes_r\Gamma_1\to\mathcal M_2\), with \(\Psi(u_g)=v_g\).
If \(g\neq e\), freeness gives \(q_{g,e}=0\); therefore \(\mathbb E_2(v_g)=0\).
On Fourier polynomials it follows that \(\tau_{H\mu}\Psi=\tau_\mu\), and continuity extends this identity to the crossed product.
Thus \(\Psi\) is a uniform \(2\)-isometry and extends, with its operator-norm bounds, to \(\mathcal M_1\).

Construct \(\Xi:\mathcal M_2\to\mathcal M_1\) from \(\theta^{-1}\) in the same way, and write
\[
 \Xi(u_h)=\sum_k u_k r_{h,k},\qquad
 r_{h,k}=\mathbf1_{\{x:\theta^{-1}(h\theta(x))=kx\}},\qquad
 p_{g,h}=\Phi_D^{-1}(q_{g,h}).
\]
On \(p_{g,h}\), one has \(\theta^{-1}(h\theta(x))=gx\).
Freeness implies \(r_{h,k}p_{g,h}=0\) for \(k\neq g\) and \(r_{h,g}p_{g,h}=p_{g,h}\).
The bounded uniform \(2\)-sums consequently give
\[
 \Xi(v_g)=\sum_h\Xi(u_h)p_{g,h}
 =\sum_hu_gp_{g,h}=u_g.
\]
The two maps are inverse on the diagonals, and the symmetric calculation handles the group unitaries of \(\mathcal M_2\).
Uniform \(2\)-density now makes them inverse on the whole completions.
\end{proof}

\begin{rem}
The preceding equivalence requires neither SBP nor comparison.
The small boundary property enters Theorem~\ref{thm:AF-model}; no comparison assumption is used there.
\end{rem}

\begin{cor}
The actions in Theorem~\ref{thm:main} are uniform measure orbit
 equivalent for every affine homeomorphism of their invariant-measure
 simplexes.
\end{cor}
\begin{proof}
Combine Theorems~\ref{thm:main} and~\ref{thm:uniform-measure-oe}.
Together with Corollary~\ref{cor:uniform-oe-full-group}, this also
proves Corollary~\ref{cor:intro-uniform-measure-oe}.
\end{proof}

\begin{rem}[Reconstruction from the uniform full group]
\label{rem:uniform-full-group-reconstruction-action}
For the actions in Theorem~\ref{thm:main},
Corollary~\ref{cor:uniform-oe-full-group} gives an isometric
full-group isomorphism from every uniform measure orbit equivalence.
Conversely, by~\cite[Theorem~8.17]{UR}, every topological group isomorphism between
these uniform full groups is implemented by a uniform measure orbit
equivalence, uniquely modulo equality almost everywhere for every
invariant measure, and is automatically isometric.
Together with the classification above, this makes affine
homeomorphism of the invariant-measure simplexes, uniform measure
orbit equivalence, and topological or isometric isomorphism of the
uniform full groups equivalent.
\end{rem}

\begin{rem}
With a unique invariant measure the bounded completion is the usual
tracial von Neumann algebra. For the present infinite free minimal
amenable systems it is the hyperfinite $\mathrm{II}_1$ factor, and
its completed diagonal is $L^\infty(X,\mu)$. Thus the result agrees
with the measured Cartan uniqueness theorem of
Connes--Feldman--Weiss~\cite{CFW}. In the multiple-trace case, the
observable and cocycle-projection requirements retain the uniform
tracial topology simultaneously over the whole simplex.
\end{rem}

\subsection{A uniform Connes--Feldman--Weiss consequence}

\begin{cor}[A Cantor minimal model and a single orbit generator]
\label{cor:uniform-CFW-action}
Let $\Gamma\curvearrowright X$ satisfy the hypotheses of
Theorem~\ref{thm:main}, and put $K=M_\Gamma(X)$.
There is a Cantor minimal homeomorphism $V:Y\to Y$ such that
$K\cong_{\mathrm{aff}}M_V(Y)$. Every affine homeomorphism
\[
 H:K\longrightarrow M_V(Y)
\]
is implemented by a uniform measure orbit equivalence
$\theta:X^0\to Y^0$ between invariant universally conull Borel sets.
In particular, there is $T\in[R_\Gamma]^u_K$ such that
\[
 \{T^nx:n\in\mathbb Z\}=\Gamma x\qquad(x\in X^0),
 \qquad T=\theta^{-1}V\theta\quad\text{on }X^0.
\]
The completed Cartan pair of the original action is consequently
isomorphic to
\[
 \bigl(\mathrm{C}(Y)^{M_V(Y)},(\mathrm{C}(Y)\rtimes_{r,V}\mathbb Z)^{M_V(Y)}\bigr),
\]
with the prescribed trace correspondence $H$.
\end{cor}

\begin{proof}
Theorem~\ref{thm:AF-model} supplies a simple unital infinite-dimensional
AF algebra $B$, its standard diagonal $E=\mathrm{C}(Y)$, and an affine
homeomorphism $K\to\T(B)$. View $Y$ as the path space of its Bratteli
diagram. Telescope the diagram so that every incidence multiplicity
is at least two. Give it a proper ordering, with unique minimal and
maximal paths. Its Vershik map $V$ is a Cantor minimal homeomorphism;
see \cite[Proposition~1.1]{GPSA}.

For completeness, its invariant measures are precisely the central
path measures corresponding to $\T(B)$. Write
$B_n=\bigoplus_v\mathrm{M}_{k_{n,v}}$ and let $c_n$ cyclically permute the
ordered length-$n$ paths ending at each vertex. For a central measure
$\nu$, all prefix cylinders ending at $v$ have equal mass. Since $c_n$
and $V$ agree except on the maximal prefix in each block,
\[
 \nu\{c_n\ne V\}
 \leq\sum_v\frac{\nu(z_{n,v})}{k_{n,v}}
 \leq\frac{1}{\min_v k_{n,v}}\longrightarrow0,
\]
where $z_{n,v}$ is the union of the prefix cylinders ending at $v$.
Every $c_n$ preserves $\nu$, so the displayed estimate, applied to
bounded continuous test functions, shows that $V$ preserves $\nu$.
Conversely, a substitution between two prefixes ending at the same
vertex is a restriction of a power of $V$: the exponent is the
difference of their ranks in the ordered prefix list. A $V$-invariant
measure therefore gives these cylinders equal mass and is central.
This identifies $M_V(Y)$ affinely homeomorphically with $\T(B)$.

The $\mathbb Z$-action generated by $V$ is free and minimal, and has
the small boundary property because $Y$ is Cantor. Apply
Theorem~\ref{thm:main} to the original action and this action, with
the specified affine homeomorphism $H$. The resulting marked
completed-pair isomorphism gives the asserted uniform measure orbit
equivalence by Theorem~\ref{thm:uniform-measure-oe}.
The sets $X^0,Y^0$ can be taken invariant, as in that theorem.
Conjugating $V|_{Y^0}$ by $\theta$ gives $T$; extend it by the identity
on $X\setminus X^0$. Since $\theta$ sends entire $\Gamma$-orbits onto
entire $V$-orbits, the displayed orbit equality holds on $X^0$.
Finally Corollary~\ref{cor:uniform-oe-full-group} transports the uniform
measure full groups. As $V$ belongs to the full group of its own
action, its transport $T$ belongs to $[R_\Gamma]^u_K$.
\end{proof}

For these actions, with $K=M_\Gamma(X)$, the uniform full group is a L\'evy group and is
extremely amenable \cite[Theorem~8.12]{UR}. Its topology determines
the relation and the affine marking of the full invariant simplex
\cite[Theorem~8.17]{UR}.

\section{Approximate conjugacy of diagonal interval maps}\label{sec:interval-maps}

In this section write $A^u$ and $D^u$ for the bounded uniform tracial
completions.  We use the following hypothesis:
\begin{equation}\label{eq:AF-pair-hypothesis}
 (D^u,A^u)\cong(E^u,B^u)
\end{equation}
by a trace-preserving isomorphism, where $B$ is a unital AF algebra and
$E$ is the diagonal of a regular unital Bratteli presentation.  In
particular, every permutation unitary at a finite stage normalizes $E$
and $E^u$.  The dynamical pairs covered by the AF model theorem satisfy
\eqref{eq:AF-pair-hypothesis}.  No lower bound on the spectral measures
is needed for the result below. We write $\mathcal U_0(A)$ for the
norm-connected component of the identity in the unitary group of $A$.

\begin{lem}\label{lem:late-block-traces}
Let $B=\varinjlim(B_n,\iota_n)$ be a unital AF algebra.  If
$x_1,\ldots,x_r\in(B_n)_{\mathrm{sa}}$ and
\[
 \sup_{\tau\in \T(B)}|\tau(x_j)|<\delta\quad(1\leq j\leq r),
\]
then, for all sufficiently large $m$, every tracial state $\sigma$ of
$B_m$ satisfies
\[
 |\sigma(\iota_{n,m}(x_j))|<\delta\quad(1\leq j\leq r).
\]
\end{lem}

\begin{proof}
The sets $K_m:=\iota_{n,m}^*(\T(B_m))\subseteq \T(B_n)$ are compact and
decreasing, and their intersection is precisely the restriction of
$\T(B)$ to $B_n$.  The latter assertion follows from compactness of the
inverse limit of the finite-stage trace simplexes, or directly by a
diagonal subsequence of compatible restrictions.  The open set
$\{\sigma:|\sigma(x_j)|<\delta\text{ for all }j\}$ contains the
intersection; compactness therefore implies that it contains $K_m$
for all sufficiently large $m$.
\end{proof}

\begin{lem}\label{lem:diagonal-matching}
For every $\eta>0$ there are a finite set
$\mathcal H\subset \mathrm{C}([0,1])$ of nonnegative contractions and $\delta>0$
with the following property.  If $a,b\in(E^u)^+_1$ satisfy
\[
 |\tau(h(a))-\tau(h(b))|<\delta
 \quad(h\in\mathcal H,\ \tau\in \T(B)),
\]
then a permutation unitary $w$ at some finite stage of $B$ satisfies
\[
 \|a-wbw^*\|_{2,\T(B)}<\eta.
\]
\end{lem}

\begin{proof}
Choose $N$ so large that $2/N<\eta^2/16$.  Put $t_j=j/N$ for
$0\leq j\leq N$ and let $h_j$ be the piecewise linear hat functions
on this mesh: $h_j(t_k)=\delta_{jk}$ and
$\sum_{j=0}^Nh_j=1$.  Thus $h_j$ is $N$-Lipschitz and
$h_j(t)\neq0$ implies $|t-t_j|\leq1/N$.
Take $\mathcal H=\{h_0,\ldots,h_N\}$ and choose $\delta>0$ so that
\[
 (N+1)\delta<\eta^2/16.
\]
Choose positive contractions $a_0,b_0$ in the diagonal of one finite
stage $B_n$ with
\[
 \|a-a_0\|_{2,\T(B)},\ \|b-b_0\|_{2,\T(B)}<\xi,
\]
where $\xi<\eta/8$ and $2N\xi<\delta/4$.  Since all these elements
belong to the same commutative algebra,
\[
 |\tau(h_j(a_0)-h_j(b_0))|<\delta+2N\xi<5\delta/4.
\]
By Lemma~\ref{lem:late-block-traces}, at a sufficiently late stage
$B_m=\bigoplus_s\mathrm{M}_{q_s}$ every normalized block trace satisfies
\[
 |\operatorname{tr}_{q_s}(h_j(a_0)-h_j(b_0))|<2\delta.
\]
Here and below $a_0,b_0$ also denote their images in $B_m$.

Fix a block and let $\mu,\nu$ be the empirical distributions of the
diagonal entries of $a_0,b_0$ in this block.  The probability measure
$\sum_j\mu(h_j)\delta_{t_j}$ can be coupled with $\mu$ so that all
transport distances are at most $1/N$; the same holds for $\nu$.
The total variation distance between these two mesh measures is
$\frac12\sum_j|\mu(h_j)-\nu(h_j)|$.  Hence their optimal mean
transport distance obeys
\[
 W_1(\mu,\nu)\leq\frac2N+
 \frac12\sum_{j=0}^N|\mu(h_j)-\nu(h_j)|
 <\frac2N+(N+1)\delta<\eta^2/8.
\]
List the entries of each matrix in increasing order and match them
in that order. If $a_1\leq\cdots\leq a_q$ and $b_1\leq\cdots\leq b_q$
are the two lists, their ordered matching has mean cost
\[
 \frac1q\sum_{i=1}^q |a_i-b_i|
 =\int_0^1|F_\mu(t)-F_\nu(t)|\,dt,
\]
where $F_\mu,F_\nu$ are the distribution functions. Indeed, at each
$t$ the two sets of indices with entry at most $t$ are nested initial
segments. For any coupling of $\mu$ and $\nu$, the probability that
its two entries lie on opposite sides of $t$ is at least
$|F_\mu(t)-F_\nu(t)|$. Integrating shows that the ordered matching
attains $W_1(\mu,\nu)$.  Since all entries belong to $[0,1]$, their mean square
distance is at most their mean absolute distance.  Thus a
permutation matrix $w_s$ in this block satisfies
\[
 \|a_0-w_sb_0w_s^*\|_{2,\operatorname{tr}_{q_s}}^2<\eta^2/8.
\]
Set $w=\bigoplus_sw_s$.  Taking the supremum over traces of $B$ and
adding the two approximation errors gives
\[
 \|a-wbw^*\|_{2,\T(B)}<2\xi+\eta/\sqrt8<\eta.
\]
By regularity of the Bratteli presentation, $w$ normalizes $E^u$.
\end{proof}

\Needspace{17\baselineskip}
\begin{thm}\label{thm:interval-uniqueness}
Assume \eqref{eq:AF-pair-hypothesis}.  For every finite set
$\mathcal F\subset \mathrm{C}([0,1])$ and $\varepsilon>0$ there are a finite
set $\mathcal H\subset \mathrm{C}([0,1])^+_1$ and $\delta>0$ such that the
following holds.  If $\phi,\psi:\mathrm{C}([0,1])\to D$ are unital
homomorphisms and
\[
 |\tau(\phi(h))-\tau(\psi(h))|<\delta
 \quad(h\in\mathcal H,\ \tau\in \T(A)),
\]
then there is a unitary $u\in\mathcal U_0(A)$ such that
\[
 \|\phi(f)-u\psi(f)u^*\|_{2,\T(A)}<\varepsilon
 \quad(f\in\mathcal F),
\]
\[
 \sup_{d\in D_1}\operatorname{dist}_{2,\T(A)}(udu^*,D_1)
 <\varepsilon,\qquad
 \sup_{d\in D_1}\operatorname{dist}_{2,\T(A)}(u^*du,D_1)
 <\varepsilon.
\]
In fact, $u$ can be chosen as $u=\exp(ic)$ with
$c\in A_{\mathrm{sa}}$ and $\|c\|\leq\pi$.
\end{thm}

\begin{proof}
Let $a=\phi(\operatorname{id}_{[0,1]})$ and
$b=\psi(\operatorname{id}_{[0,1]})$.  Through the isomorphism in
\eqref{eq:AF-pair-hypothesis}, view $a,b$ as positive contractions
in $E^u$.  Choose $\eta>0$ so small that whenever $x,y$ are
commuting positive contractions with $\|x-y\|_{2,\T(B)}<\eta$,
\[
 \|f(x)-f(y)\|_{2,\T(B)}<\varepsilon/3\quad(f\in\mathcal F).
\]
This follows directly from uniform continuity: for any $r>0$,
\[
 \|f(x)-f(y)\|_{2,\T(B)}
 \leq\omega_f(r)+2\|f\|\|x-y\|_{2,\T(B)}/r.
\]
Use the tests and tolerance supplied by
Lemma~\ref{lem:diagonal-matching}.  It yields a finite-stage permutation
unitary $w\in B_m$ for which $\|a-wbw^*\|_{2,\T(B)}<\eta$; note that
$a$ and $wbw^*$ commute.  Pulling $w$ back gives a unitary $v\in A^u$
normalizing $D^u$, with
\[
 \|\phi(f)-v\psi(f)v^*\|_{2,\T(A)}<\varepsilon/3
 \quad(f\in\mathcal F).
\]
Since $B_m$ is finite dimensional, spectral calculus gives
$h\in(B_m)_{\mathrm{sa}}$ with $\|h\|\leq\pi$ and $w=\exp(ih)$.
Let $k\in(A^u)_{\mathrm{sa}}$ be its pullback. Then
$\|k\|\leq\pi$ and $v=\exp(ik)$.
Choose $\zeta>0$ with
\[
 2\zeta<\varepsilon,\qquad
 2\zeta\max\{1,\|f\|:f\in\mathcal F\}<\varepsilon/3.
\]
Bounded density, followed by symmetrization, gives
$c\in A_{\mathrm{sa}}$ with $\|c\|\leq\pi$ and
$\|c-k\|_{2,\T(A)}<\zeta$. Put $u=\exp(ic)$.
This is a unitary in $\mathcal U_0(A)$, with path
$t\mapsto\exp(itc)$ from $1$ to $u$.
The Duhamel identity
\[
 \exp(ic)-\exp(ik)
 =i\int_0^1\exp(i(1-t)c)(c-k)\exp(itk)\,dt
\]
and invariance of each tracial $2$-norm under multiplication by
unitaries give
\[
 \|u-v\|_{2,\T(A)}\leq\|c-k\|_{2,\T(A)}<\zeta.
\]
For every contraction $d\in D$,
\[
 \|udu^*-vdv^*\|_{2,\T(A)}\leq2\zeta,
 \qquad
 \|u^*du-v^*dv\|_{2,\T(A)}\leq2\zeta.
\]
The two elements on the right belong to $(D^u)_1$, which is the
uniform $2$-norm closure of $D_1$.  These bounds therefore prove
both uniform distance assertions. For $f\in\mathcal F$,
\[
 \|u\psi(f)u^*-v\psi(f)v^*\|_{2,\T(A)}\leq2\zeta\|f\|.
\]
Combining the estimates proves the theorem.
\end{proof}

\begin{rem}
If the trace spectral distributions of $\phi$ and $\psi$ agree
exactly, the theorem can be applied for every tolerance.  It gives
approximate conjugacy by unitaries in $\mathcal U_0(A)$ that uniformly
almost normalize the whole diagonal unit ball in the uniform $2$-norm.
The logarithm argument
uses the finite-stage permutation supplied by
Lemma~\ref{lem:diagonal-matching}, and requires no stable-rank
assumption on $A$. It does not assert that a single unitary in $A$ exactly
conjugates the two maps.
\end{rem}

\section{Boolean algebras of the completed diagonals}
\label{sec:boolean-completions}

Let $X$ be compact metrizable and let $K\subseteq\Prob(X)$ be
nonempty and weak* compact. We use the bounded completion
$D_K=\mathrm{C}(X)^K$, with the common null kernel removed, and put
\[
 \B_K=\Proj(D_K),\qquad
 c_K(p)=\sup_{\mu\in K}\mu(p),\qquad
 d_K(p,q)=c_K(p\symdiff q)=\|p-q\|_{2,K}^{\,2}.
\]
Here $p\symdiff q=p(1-q)+(1-p)q$. Finite Boolean operations
are continuous for $d_K$. We distinguish the bare Boolean algebra,
its exact capacity or disagreement metric, and the marked
completed diagonal $(D_K,K)$. Metric completeness concerns
$d_K$; order completeness concerns joins. A Boolean algebra is
$\sigma$-complete if every countable subset has a supremum, and
complete if every subset does. Only the statements explicitly
restricted to full-trace AF models below use the classification
theorems of this paper.

\subsection{Uniform completion and norm real rank zero}\label{subsec:boolean-RR0}

Elements of $D_K$ have bounded Borel representatives modulo common
$K$-null sets. A bounded uniform $2$-Cauchy sequence has a subsequence
with summable successive uniform $2$-differences. That subsequence
converges pointwise outside one Borel set null for every $\mu\in K$,
and its limit is the required representative. Projections are
characteristic functions modulo common null sets.

For $a\in D_K$, the functions $\mu\mapsto\mu(a)$ and
$\mu\mapsto\mu(a^*a)$ are continuous on $K$, by bounded uniform
approximation. Multiplication is continuous on bounded sets.
Thus $\B_K$ is closed in the uniformly complete unit ball and is
complete for $d_K$. It is separable because that unit ball is the
bounded uniform completion of the norm-separable continuous diagonal.

\begin{lem}\label{lem:boolean-threshold}
Suppose increasing finite-dimensional commutative algebras
$F_n\subseteq D_K$ have boundedly uniform $2$-dense union.
Their projection algebras have $d_K$-dense union in $\B_K$.
\end{lem}
\begin{proof}
For $p\in\B_K$, choose positive step contractions $h_n\in F_n$ with
$\norm{h_n-p}_{2,K}\to0$. The projections
$e_n=1_{[1/2,1]}(h_n)$ satisfy the pointwise inequality
$|e_n-p|\leq2|h_n-p|$, and therefore converge to $p$ uniformly in
$2$-norm.
\end{proof}

Relative $K$-SBP means that $X$ has a basis of open sets whose
boundaries are null for all $\mu\in K$.

\begin{lem}[Finite diagonal approximation]\label{lem:boolean-finite-diagonals}
If $(X,K)$ has $K$-SBP, there are increasing finite-dimensional
unital subalgebras $F_n\subseteq D_K$ whose union is boundedly
uniformly $2$-dense in $D_K$.
\end{lem}
\begin{proof}
For a small-boundary open set $\varnothing\ne U\ne X$, put
\[
 a_n(x)=\min\{1,n\,d(x,X\setminus U)\},\qquad
 b_n(x)=\max\{0,1-n\,d(x,\overline U)\}.
\]
Then $a_n\leq1_U\leq b_n$, and $b_n-a_n$ decreases to
$1_{\partial U}$. The continuous functions
$\mu\mapsto\int(b_n-a_n)\,d\mu$ decrease to zero on compact $K$.
Dini's theorem gives $\|1_U-a_n\|_{2,K}\to0$, so $1_U\in D_K$.
Take finite covers by such sets with mesh tending to zero, turn
them into Borel partitions by ordered differences, and successively
refine the partitions. Their characteristic functions generate
nested finite diagonals $F_n$. Sampling a continuous contraction
on their pieces gives step contractions converging in norm.
The unit-ball description of the bounded completion then gives
the required bounded uniform $2$-density.
\end{proof}

\begin{lem}[From uniform approximation to real rank zero]
\label{lem:boolean-RR0}
Let $D$ be a commutative unital $\mathrm{C}^*$-algebra and let $K$ be a nonempty
family of states such that
\[
 \|a\|_{2,K}=\sup_{\tau\in K}\tau(a^*a)^{1/2}
\]
is a norm and the closed unit ball of $D$ is complete for this norm.
If the finite-spectrum self-adjoint contractions are uniformly
$2$-norm dense in $(D_{\mathrm{sa}})_1$, then $D$ has real rank zero
in its $\mathrm{C}^*$-norm.  In particular, finite-spectrum self-adjoint
elements are norm dense in $D_{\mathrm{sa}}$.
\end{lem}

\begin{proof}
Fix $a=a^*\in D$ and $\eta>0$.  We may assume $a\ne0$.
By the hypothesis and scaling, choose finite-spectrum self-adjoint
elements $b_n$ such that
\[
 \|b_n\|\leq\|a\|,\qquad
 \|b_n-a\|_{2,K}\leq\eta2^{-3n}\qquad(n\geq1).
\]
Put
\[
 t_n=\frac{\eta}{2}(1-2^{-n}),\qquad
 q_n^+=1_{[t_n,\infty)}(b_n),\qquad
 q_n^-=1_{(-\infty,-t_n]}(b_n).
\]
These spectral projections belong to $D$ because $b_n$ has finite
spectrum.  Define recursively
\begin{equation}\label{eq:boolean-RR0-projection-recursion}
 p_1=q_1^+,\qquad
 p_n=q_n^++(1-q_n^+-q_n^-)p_{n-1}\quad(n\geq2).
\end{equation}
All elements commute, so every $p_n$ is a projection and
\[
 q_n^+\leq p_n\leq1-q_n^-.
\]
The recursion retains the previous projection inside the middle
spectral interval and prescribes its value outside that interval.

At a character where $p_{n-1}=0$ and $p_n=1$, one has
$b_{n-1}<t_{n-1}$ and $b_n\geq t_n$.  At a character where
$p_{n-1}=1$ and $p_n=0$, one has $b_{n-1}>-t_{n-1}$ and
$b_n\leq-t_n$.  Consequently, pointwise on the spectrum of $D$,
\[
 (t_n-t_{n-1})|p_n-p_{n-1}|\leq|b_n-b_{n-1}|.
\]
It follows that
\begin{align*}
 \|p_n-p_{n-1}\|_{2,K}
 &\leq\frac{\|b_n-b_{n-1}\|_{2,K}}{t_n-t_{n-1}}\\
 &\leq\frac{9\eta2^{-3n}}{\eta2^{-n-1}}
 =18\,2^{-2n}\qquad(n\geq2).
\end{align*}
Thus $(p_n)$ is uniformly $2$-norm Cauchy.  Completeness and
continuity of multiplication on bounded sets give a projection
$p\in D$ with $\|p_n-p\|_{2,K}\to0$.

The defining inequalities also give
\[
 (b_n-t_n)_+(1-p_n)=0,\qquad
 (-b_n-t_n)_+p_n=0.
\]
Since $D$ is commutative, the positive-part function is
$1$-Lipschitz for $\|\cdot\|_{2,K}$ on self-adjoint elements.
Passing to the limit, and using faithfulness of this norm, yields
\[
 (a-\eta/2)_+(1-p)=0,\qquad
 (-a-\eta/2)_+p=0.
\]
Hence $ap\geq-(\eta/2)p$ and
$a(1-p)\leq(\eta/2)(1-p)$.  The self-adjoint element
\[
 c=a+\eta(2p-1)
\]
therefore satisfies
\[
 cp\geq(\eta/2)p,\qquad
 c(1-p)\leq-(\eta/2)(1-p).
\]
In particular, $|c|\geq(\eta/2)1$, so $c$ is invertible, and
$\|c-a\|=\eta$.  Taking $\eta$ smaller than any prescribed
tolerance proves norm density of the invertible self-adjoint
elements, which is real rank zero.
\end{proof}

\begin{thm}[Small boundaries and real rank zero]
\label{thm:boolean-SBP-RR0}
Let $X$ be compact metrizable and let $K\subseteq\Prob(X)$ be
nonempty and compact.  Then
\[
 (X,K)\text{ has }K\text{-SBP}
 \quad\Longleftrightarrow\quad
 D_K=\mathrm{C}(X)^K\text{ has real rank zero}.
\]
Here real rank zero refers to the $\mathrm{C}^*$-norm of the bounded uniform
completion, after quotienting by the common trace kernel.
\end{thm}

\begin{proof}
Under $K$-SBP, Lemma~\ref{lem:boolean-finite-diagonals}
supplies uniformly $2$-norm dense finite-spectrum self-adjoint
contractions.  The bounded uniform completion has faithful
$\|\cdot\|_{2,K}$ and complete closed unit ball, so
Lemma~\ref{lem:boolean-RR0} applies.
Conversely, norm density of invertible self-adjoint elements gives
their uniform $2$-norm density, and
\cite[Theorem~2.16 and Corollary~2.17]{EN} gives $K$-SBP.
\end{proof}

\begin{rem}[The real rank of the entire completion]
\label{rem:boolean-completion-RR0-question}
Theorem~2.16 of~\cite{EN} proves that, under $K$-SBP, self-adjoint
elements of the original continuous diagonal $\mathrm{C}(X)$ admit norm
approximation by invertible self-adjoint elements of $\mathrm{C}(X)^K$.
Remark~2.18 of that paper asks whether $K$-SBP is equivalent to
real rank zero of the entire completion.
Theorem~\ref{thm:boolean-SBP-RR0} answers this question
affirmatively.  The additional step is the summable projection
construction in~\eqref{eq:boolean-RR0-projection-recursion}, which
applies to every self-adjoint element of the completion.
Neither Lemma~\ref{lem:boolean-RR0} nor
Theorem~\ref{thm:boolean-SBP-RR0} requires amenability,
aperiodicity, or a Bauer hypothesis on the designated family.
\end{rem}

We recall the following well-known form of Stone duality.

\begin{prop}[Boolean reconstruction]\label{prop:boolean-Stone}
If $D$ is commutative, unital, and has norm real rank zero, then
\[
 D\cong \mathrm{C}(\Stone(\Proj D)).
\]
Every Boolean isomorphism of projection algebras extends uniquely to a
unital norm $\mathrm{C}^*$-isomorphism. Every finitely additive
Boolean probability extends uniquely to a state of $D$.
\end{prop}
\begin{proof}
Finite linear combinations of projections are norm dense. A Boolean
isomorphism preserves finite partitions, so its linear extension
preserves the maximum absolute coefficient over nonzero partition
pieces, which is the norm. It therefore extends uniquely. Boolean
probabilities similarly extend to positive normalized functionals.
Stone duality identifies the resulting algebra
\cite[Chapters~31 and~36]{Fremlin}.
\end{proof}

The Stone space here is generally not metrizable. Uniform
$2$-separability does not assert norm separability.

\subsection{An intrinsic range theorem}\label{subsec:boolean-range}

\begin{defn}
A \emph{normalized upper probability} on $\B$ is a function
$c:\B\to[0,1]$ of the form $c(p)=\sup_{\nu\in\mathcal P}\nu(p)$,
where $\mathcal P$ is a nonempty family of finitely additive
probabilities. It is \emph{strictly positive} if $c(p)>0$ for $p\ne0$.
Its \emph{dominated probability core} is
\[
 K_c=\{\nu\in\Probfa(\B):\nu(p)\leq c(p)\text{ for every }p\in\B\}.
\]
\end{defn}

An upper probability is monotone and subadditive, with $c(0)=0$ and
$c(1)=1$. Strict positivity makes $d_c(p,q)=c(p\symdiff q)$ a metric.
The upper-envelope hypothesis is stronger than subadditivity alone.

\begin{lem}\label{lem:boolean-core}
The set $K_c$ is nonempty compact convex in the pointwise topology,
and
\[
 c(p)=\max_{\nu\in K_c}\nu(p),\qquad
 |\nu(p)-\nu(q)|\leq d_c(p,q).
\]
If $(\B,d_c)$ is separable, restriction to a countable dense Boolean
subalgebra is injective on $K_c$, and $K_c$ is metrizable.
\end{lem}
\begin{proof}
Finite additivity, positivity, normalization, and domination are
closed conditions in $[0,1]^\B$. The original representing family
is contained in $K_c$, proving equality of upper envelopes.
Compactness gives attainment. The difference estimate follows from
finite additivity. Dense-core values determine all values, so the
continuous injection into countably many core coordinates is a
homeomorphism onto its image.
\end{proof}

\begin{thm}[Range of the projection Boolean algebras]\label{thm:boolean-range}
For a nontrivial Boolean algebra $\B$, the following are equivalent:
\begin{enumerate}
\item $\B\cong\Proj(\mathrm{C}(X)^K)$ for compact metrizable $X$ and nonempty
compact $K\subseteq\Prob(X)$;
\item $\B$ admits a normalized strictly positive upper probability $c$
for which $d_c$ is complete and separable.
\end{enumerate}
The realization preserves $c$ and can be chosen with $X$
zero-dimensional. Its completed diagonal has norm real rank zero.
\end{thm}
\begin{proof}
For (i)$\Rightarrow$(ii), use $c(p)=\sup_{\mu\in K}\mu(p)$.
Faithfulness after the common-null quotient gives strict positivity;
Subsection~\ref{subsec:boolean-RR0} gives metric completeness and separability. No SBP
assumption on the original $X$ is needed for this direction.

Conversely, generate a countable Boolean subalgebra $\B_0$ from a
countable $d_c$-dense subset of $\B$. Put $X=\Stone(\B_0)$.
Each Boolean probability on $\B_0$ extends uniquely to a Radon
probability on $X$, since clopen simple functions are norm dense in
$\mathrm{C}(X)$. Let $K$ be the restrictions of $K_c$. Pointwise convergence
on clopens is weak* convergence of these measures, so
Lemma~\ref{lem:boolean-core} makes $K$ compact metrizable. Its upper
envelope on $\B_0$ is $c$.

Clopen simple functions are boundedly uniform $2$-dense in $\mathrm{C}(X)^K$.
Lemma~\ref{lem:boolean-threshold} therefore makes $\B_0$ dense in its
projection algebra. Both that algebra and $\B$ are metric
completions of the same core. The canonical isometry is Boolean
because finite Boolean operations are continuous for $d_c$.
Finally, $X$ is zero-dimensional and hence has $K$-SBP, so
Theorem~\ref{thm:boolean-SBP-RR0} applies.
\end{proof}

\begin{rem}[Scope]
This is an existence criterion for a bare algebra, expressed by the
existence of a suitable capacity. If an original realization lacks SBP,
the construction realizes its Boolean part anew and need not recover
its original norm algebra. For a general designated family, $K_c$
need not equal that family and need not be a Choquet simplex.
Recovery of the original simplex is an additional conclusion in
Subsection~\ref{subsec:boolean-simplex}.
\end{rem}

\begin{prop}[A necessary chain condition]\label{prop:boolean-ccc}
Every Boolean algebra in Theorem~\ref{thm:boolean-range} carries a strictly
positive finitely additive probability and satisfies ccc. Its cardinality
is at most the continuum.
\end{prop}
\begin{proof}
Choose a dense sequence $(\nu_n)$ in $K_c$ and set
$\nu_0=\sum_{n\geq1}2^{-n}\nu_n$. If $\nu_0(p)=0$, continuity of
evaluation on $K_c$ gives $c(p)=0$, hence $p=0$.
Each member of a disjoint nonzero family has positive $\nu_0$-mass;
for each $n$ only finitely many can have mass at least $1/n$.
Thus the family is countable. The cardinality bound follows from
metric separability.
\end{proof}

\begin{example}[An excluded atomless algebra]\label{ex:boolean-PN}
The algebra $\mathcal P(\N)/\mathrm{Fin}$ cannot occur in
Theorem~\ref{thm:boolean-range}. Identify $\N$ with the finite binary strings.
For each infinite binary sequence, its set of prefixes is infinite;
two different branches have only finitely many common prefixes.
Their classes give an uncountable disjoint family of nonzero elements.

The algebra is atomless and not $\sigma$-complete. For the latter,
partition $\N$ into infinite sets $A_n$. From any upper bound $[B]$
of their classes, remove one point of $B\cap A_n$ for every $n$.
The result is still an upper bound and is strictly smaller modulo
finite sets.
\end{example}

Without the compact-metrizable restriction, every Boolean algebra has
a trivial action realization. On $S=\Stone(\B)$ take the trivial
$\mathbb Z$-action and all probabilities. Point masses give
$\norm{f}_{2,\Prob(S)}=\norm{f}_\infty$, so the completion is $\mathrm{C}(S)$
and its projection algebra is $\B$.

\subsection{Recovering and classifying the trace simplex}\label{subsec:boolean-simplex}

For the statements involving simplex classification, a
\emph{full-trace AF model} means a marked completed diagonal
$(D,K)$ isomorphic to $(E^u,\T(B))$, where $B$ is separable,
simple, unital, and infinite-dimensional AF and $E$ is its standard
diagonal. Theorem~\ref{thm:AF-model} places the actions considered
in this paper in this class, and Theorem~\ref{thm:AF-pairs} classifies
it by its full trace simplex. Each such model also has a free minimal
Cantor $\mathbb Z$-model: realize its metrizable Choquet simplex
by \cite{Dow}, apply Theorem~\ref{thm:AF-model}, and then
Theorem~\ref{thm:AF-pairs}.
The preceding range theorem and real rank zero lemma do not
require these additional assumptions.

Theorem~\ref{thm:AF-pairs} permits replacing a standard AF model
$A$ by $A\otimes \mathrm{M}_{2^\infty}$ with its product diagonal, since its
full trace simplex is unchanged.

\begin{thm}[Recovery from the scalar capacity]\label{thm:boolean-recovery}
For a diagonal in the full-trace AF class,
\[
 K=\{\nu\in\Probfa(\B_K):\nu(p)\leq c_K(p)
                         \text{ for every }p\in\B_K\}.
\]
This identification is affine and homeomorphic for pointwise
convergence on projections.
\end{thm}
\begin{proof}
Work in a CAR-stabilized standard AF model and let $\B_0$ be the
union of its finite diagonal projection algebras. It is $d_K$-dense.
A dominated Boolean probability $\nu$ restricts to a state of the
norm diagonal core.

For a dyadic positive step contraction in one finite core write
\[
 f=\sum_i\frac{m_i}{q}p_i,\qquad q=2^r,\quad 0\leq m_i\leq q,
\]
where the $p_i$ partition the unit. Choose equal-trace projections
$e_0,\ldots,e_{q-1}$ in $r$ later CAR coordinates. Every trace on
the tensor product is a product with the CAR trace, so
$\tau(p_i e_j)=\tau(p_i)/q$. Define, modulo $q$,
\[
 r_k=\sum_i p_i\sum_{j=0}^{m_i-1}e_{j+k}.
\]
These are core projections with
\[
 \frac1q\sum_{k=0}^{q-1}r_k=f,\qquad
 \tau(r_k)=\tau(f)\quad(\tau\in K).
\]
Consequently,
\[
 \nu(f)=\frac1q\sum_k\nu(r_k)
       \leq\sup_{\tau\in K}\tau(f).
\]
Dyadic approximation, scaling, and translation give this inequality
for every self-adjoint core element. Compact convex separation
therefore makes its restriction the restriction of some $\mu\in K$.

For arbitrary $p\in\B_K$, take $p_n\in\B_0$ with $d_K(p_n,p)\to0$.
Both probabilities satisfy the domination estimate
$|\nu(p)-\nu(p_n)|\leq c_K(p\symdiff p_n)$.
Hence equality on the core extends to every projection.
The reverse inclusion is immediate. Evaluation continuity and
compactness give the topology assertion.
\end{proof}

\begin{thm}[Metric Boolean classification]\label{thm:boolean-metric-classification}
In the full-trace AF class,
\[
 (\B_K,d_K)\cong_{\mathrm{isometric\ Boolean}}(\B_L,d_L)
 \quad\Longleftrightarrow\quad K\cong_{\mathrm{aff}}L.
\]
Every affine homeomorphism lifts to such a Boolean isomorphism.
Every Boolean isomorphism preserving the capacity induces an affine
homeomorphism of the recovered simplexes.
\end{thm}
\begin{proof}
The forward implication follows from Theorem~\ref{thm:boolean-recovery},
because $c_K(p)=d_K(p,0)$ and Boolean coordinate isomorphisms
transport dominated probability cores. The reverse implication is
Theorem~\ref{thm:AF-pairs}.
\end{proof}

\begin{rem}
Recovery uses the AF core and uniform projection density; it does not
require the additional norm RR0 argument. Norm RR0 additionally
reconstructs the \emph{whole norm algebra} from its bare Boolean
algebra. The recovered $K$ is generally much smaller than the full
state space of $D_K$.
\end{rem}

\subsection{Countable joins and order completeness}\label{subsec:boolean-joins}

\begin{lem}\label{lem:boolean-separable-obstruction}
Let $\B$ be $\sigma$-complete and suppose
$d_c(p,q)=c(p\symdiff q)$ is a separable metric from a strictly
positive submeasure. For every disjoint sequence $(a_n)$,
$c(a_n)\to0$.
\end{lem}
\begin{proof}
Otherwise pass to a subsequence with $c(a_n)\geq\eps>0$.
For every $E\subseteq\N$, form $b_E=\bigvee_{n\in E}a_n$.
If $E\ne F$, some $a_n$ is below $b_E\symdiff b_F$, so
$d_c(b_E,b_F)\geq\eps$. This contradicts metric separability.
\end{proof}

\begin{lem}\label{lem:boolean-disjoint-clopens}
For a minimal Cantor homeomorphism whose invariant simplex has
infinitely many extreme points, there are disjoint clopen $A_n$
and invariant $\nu_n$ with $\nu_n(A_n)>1/2$.
\end{lem}
\begin{proof}
Choose distinct ergodic $\mu_j$ with a subsequence converging to
$\lambda$. Let $a_j$ be the weight of $\mu_j$ in the ergodic
decomposition of $\lambda$. Mutual singularity gives
$\sum_j a_j\leq1$, hence $a_j\to0$. Let $E_j$ be the Borel set
of $\mu_j$-generic points for a countable family of continuous
functions determining probability measures. The ergodic theorem
and decomposition give $\mu_j(E_j)=1$ and $\lambda(E_j)=a_j$.

Choose $\delta_n>0$ with $\sum_n\delta_n<1/8$.
Inductively let $F$ be the union of the previous disjoint clopen
$A_i$, with $\lambda(A_i)<\delta_i$. Choose a late, unused $j$
with $\mu_j(F)<1/4$ and $a_j<\delta_n/2$.
A compact $H\subseteq E_j\setminus F$ can have $\mu_j(H)>1/2$.
Since $\lambda(H)\leq a_j$, it has a clopen neighborhood $A_n$
disjoint from $F$ with $\lambda(A_n)<\delta_n$. Take $\nu_n=\mu_j$.
\end{proof}

\begin{thm}[Sharp completeness criterion]\label{thm:boolean-completeness}
In the full-trace AF class,
\[
 \B_K\text{ is complete}
 \ \Longleftrightarrow\
 \B_K\text{ is }\sigma\text{-complete}
 \ \Longleftrightarrow\
 |\partial_e K|<\infty.
\]
For $m$ extreme points,
\[
 D_K\cong\prod_{j=1}^mL^\infty([0,1]).
\]
Its bare projection algebra is the standard atomless probability
algebra, independently of $m\geq1$.
\end{thm}
\begin{proof}
Use a Cantor minimal model. For finitely many ergodic measures
their average $\mu_0$ gives
\[
 \norm{a}_{2,\mu_0}^2\leq\norm{a}_{2,K}^2
 \leq m\norm{a}_{2,\mu_0}^2.
\]
Thus $D_K=L^\infty(X,\mu_0)$. Mutual singularity and atomlessness
give the displayed product. Its measure algebra is complete; the
finite disjoint union is again a standard atomless probability
space \cite[Sections~321 and~331]{Fremlin}.

For infinitely many extreme points, Lemma~\ref{lem:boolean-disjoint-clopens}
gives disjoint regular projections of capacity greater than $1/2$.
Lemma~\ref{lem:boolean-separable-obstruction} excludes $\sigma$-completeness.
\end{proof}

The theorem also distinguishes the completed diagonals as unmarked
norm $\mathrm{C}^*$-algebras. In the full-trace AF class,
\[
 D_K\cong L^\infty([0,1])
 \quad\Longleftrightarrow\quad |\partial_eK|<\infty.
\]
Indeed, every finite product in the theorem is isomorphic to
$L^\infty([0,1])$. With infinitely many extreme points,
$\Proj(D_K)$ is not $\sigma$-complete. A $\mathrm{C}^*$-isomorphism
restricts to an order isomorphism of projection Boolean algebras,
whereas the projections of a commutative von Neumann algebra
form a complete Boolean algebra. Consequently $D_K$ is not
$\mathrm{C}^*$-isomorphic to any commutative von Neumann algebra,
and in particular not to the $L^\infty$-algebra of a standard
probability space. For an action in Theorem~\ref{thm:main},
this applies precisely when there are infinitely many ergodic
invariant probability measures. Proposition~\ref{prop:boolean-Bauer}
gives explicit section-algebra models for Bauer simplexes.

The term \emph{infinite simplex} here means infinitely many extreme
points. All these Boolean algebras are metrically complete. In Stone
duality $\sigma$-completeness means a basically disconnected Stone
space \cite[Section~314]{Fremlin}; norm RR0 gives zero-dimensionality.

\begin{example}[Why fullness is essential]\label{ex:boolean-designated}
The preceding criterion does not hold for arbitrary designated
Bauer simplexes. On $[0,1]$ with Lebesgue measure, let $r_n$ be
independent Rademacher functions and put
\[
 h_\infty=1,\qquad h_n=1+2^{-n}r_n.
\]
The barycentres of these densities form a compact family $K$ of
measures affinely homeomorphic to the probability simplex of a
convergent sequence. Indeed, $h_n\to1$ in $L^\infty$, and the
coefficients of a barycentre are recovered from
$\int r_n h\,dm=2^{-n}\lambda_n$. Thus $K$ is a Bauer simplex with
infinitely many extreme points. Nevertheless, all densities lie
between $1/2$ and $3/2$. Its uniform $2$-norm is equivalent to
ordinary $L^2$, so its completion is $L^\infty$ and its Boolean
algebra is complete.
\end{example}

The bare algebra and even its uniform topology lose the finite
vertex count: all finite products in Theorem~\ref{thm:boolean-completeness}
have the topology of their single averaged measure. The exact
metric retains the simplex by Theorem~\ref{thm:boolean-metric-classification}.

\subsection{Bauer models and different bare Boolean algebras}\label{subsec:boolean-Bauer}

Let $\Omega=\{0,1\}^{\N}$ with fair product measure $m$, and let $S$
be the binary odometer. For compact metrizable $Z$ define
\[
 T_Z(z,\omega)=(z,S\omega),\qquad X_Z=Z\times\Omega.
\]
This action is free. Unique ergodicity of the odometer and
disintegration over $Z$ give
\[
 K_Z=M_{T_Z}(X_Z)
     =\{\lambda\otimes m:\lambda\in\Prob(Z)\}.
\]

\begin{prop}[Bauer model]\label{prop:boolean-Bauer}
The action has SBP and
\[
 \mathrm{C}(X_Z)^{K_Z}\cong \mathrm{C}_\sigma(Z,L^\infty(\Omega,m)),
\]
where the right side consists of uniformly norm-bounded
$L^2$-continuous sections. If $\Aprob$ is the standard atomless
probability algebra with its measure metric, then
\[
 \B_{K_Z}\cong \mathrm{C}(Z,\Aprob)
\]
with pointwise Boolean operations and uniform measure metric.
\end{prop}
\begin{proof}
The uniform $2$-norm on continuous functions is
$\sup_z\norm{f(z,\cdot)}_{2,m}$. A bounded uniform $2$-norm limit yields an
$L^2$-continuous section. Conversely, approximate such a section
locally by finitely many bounded continuous fibre functions and
combine them with a partition of unity on $Z$.

For SBP, approximate a continuous real $f$ uniformly by $f_0$
depending on $z$ and the first $N$ binary coordinates.
The continuous tail function
\[
 h(\omega)=\sum_{j\geq1}2^{-j}\omega_{N+j}
\]
is uniform on $[0,1]$, independently of those first coordinates.
For $\alpha>0$, $g=f_0+\alpha h$ satisfies
\[
 m\{\omega:|g(z,\omega)|<\delta\}\leq2\delta/\alpha
 \quad(z\in Z).
\]
Its zero set is null for every $\lambda\otimes m$.
Start with a function negative at a specified point and positive
on the complement of a prescribed neighborhood. Small such
perturbations give a small-boundary neighborhood. The projection
description follows from the section identification.
\end{proof}

The same section algebra is the marked completed diagonal of every
full-trace AF model with simplex $\Prob(Z)$. By \cite{Dow} and
Theorems~\ref{thm:AF-model} and~\ref{thm:AF-pairs}, it suffices to
prove this for a minimal Cantor homeomorphism $T:Y\to Y$ with ergodic
measures $\mu_z$ indexed homeomorphically by $z\in Z$.
Identify $(\Omega,m)$ as a probability space with the Lebesgue unit
interval. Choose refining clopen partitions of $Y$ of vanishing mesh,
ordering the children of each atom. At $z$, assign to each atom $P$
an interval of length $\mu_z(P)$, subdividing the parent's interval
in this order. The interval endpoints vary continuously with $z$.
Sending $1_P$ to the indicator of this interval defines a unital
homomorphism on the clopen step functions into
$\mathrm{C}_\sigma(Z,L^\infty[0,1])$. It preserves each fibre integral and
hence the uniform $2$-norm, since every invariant measure is a
barycentre of the $\mu_z$. It extends to the bounded completion.
At each $z$ the interval meshes tend to zero by atomlessness and
compactness of $Y$, so images of clopen step contractions are
$L^2$-dense in the fibre unit ball.

The image also contains the scalar sections $\mathrm{C}(Z)1$. Indeed, for
$g\in \mathrm{C}(Y)$ put $g_n=n^{-1}\sum_{j=0}^{n-1}g\circ T^j$ and
$h(z)=\mu_z(g)$. The continuous functions
$z\mapsto\|g_{2^k}-h(z)\|_{2,\mu_z}^2$ decrease with $k$, by
invariance and convexity, and tend pointwise to zero by ergodicity.
Dini's theorem gives uniform convergence. Thus the image contains
$h1$. The functions $z\mapsto\mu_z(g)$ separate points and
are closed under conjugation. By Stone--Weierstrass they generate
$\mathrm{C}(Z)$ as a norm-closed unital algebra. Given a section in the unit ball, choose at each point a
clopen step contraction approximating its value in $L^2$. The estimate
persists in a neighborhood by continuity. A finite partition of unity
on $Z$ combines these local approximants into one contraction in the
image with the same uniform error. Completeness on the unit ball gives
surjectivity. This proves the marked identification and applies, by
Theorem~\ref{thm:main}, to every action covered by that theorem
with invariant simplex $\Prob(Z)$. Alternatively, the product action
above has infinite orbits, relative SBP, and its full invariant family,
so the identification also follows from uniform relation
classification~\cite[Theorem~8.6 and Corollary~8.7]{UR}.

\begin{example}[An explicit missing join]\label{ex:boolean-sequence}
For $Z=\{0\}\cup\{1/n:n\geq1\}$, the diagonal is
\[
 D_Z=\{(a_0,a_1,\ldots):
       \sup_n\norm{a_n}<\infty,\
       \norm{a_n-a_0}_{2,m}\to0\}.
\]
Let $e_k$ be the scalar projection supported at $1/(2k)$.
A supremum $p$ would equal $1$ at every even coordinate,
forcing $p(0)=1$. For each odd coordinate, the projection equal
to $1$ except there is an upper bound of all $e_k$.
Leastness would force $p$ to vanish at all odd coordinates,
contradicting continuity at $0$.
\end{example}

\begin{prop}[No complete corner over a Cantor base]\label{prop:boolean-no-corner}
No nonzero principal ideal of $\mathrm{C}(C,\Aprob)$, where $C$ is Cantor,
is $\sigma$-complete.
\end{prop}
\begin{proof}
Fix $p\ne0$ and $z_0$ with $p(z_0)\ne0$. Choose strictly
decreasing clopen neighborhoods $U_n$ of $z_0$ with diameters
tending to zero and nonempty shells $V_n=U_n\setminus U_{n+1}$.
Put $a_n=p1_{V_{2n}}$.
If $q=\bigvee_na_n$ existed in $[0,p]$, it would equal $p$
on every even shell, forcing $q(z_0)=p(z_0)$.
For each odd shell, $p1_{C\setminus V_{2n-1}}$ is an upper
bound. Leastness forces $q=0$ there, hence $q(z_0)=0$ by
continuity along odd shells, a contradiction.
\end{proof}

\begin{prop}[Different bare Boolean algebras]\label{prop:boolean-different}
Let
\[
 K_0=\Prob(C),\qquad K_1=\Prob(\{*\}\sqcup C).
\]
Their completed projection algebras are not isomorphic as bare
Boolean algebras. Consequently, the completed diagonals
$D_{K_0}$ and $D_{K_1}$ are not isomorphic as norm
$\mathrm{C}^*$-algebras. Both Boolean algebras are
non-$\sigma$-complete, and both simplexes have uncountably many
extreme points.
\end{prop}
\begin{proof}
Their Bauer models are
\[
 \B_{K_0}=\mathrm{C}(C,\Aprob),\qquad \B_{K_1}\cong\Aprob\times \mathrm{C}(C,\Aprob).
\]
In the second algebra $e=(1,0)$ has $[0,e]\cong\Aprob$, which is
complete. Proposition~\ref{prop:boolean-no-corner} excludes every such
nonzero corner in the first. Both whole algebras have a
non-$\sigma$-complete principal ideal isomorphic to $\mathrm{C}(C,\Aprob)$.
Any $\mathrm{C}^*$-isomorphism of the completed diagonals would
restrict to a Boolean isomorphism of their projections, so their
norm algebras are also nonisomorphic.
\end{proof}

The product actions here are explicit free Cantor actions, though
nonminimal. Downarowicz's theorem \cite{Dow} supplies minimal
Cantor systems with these two simplexes; their SBP is automatic.
The preceding marked section-algebra identification applies to
both minimal realizations. Thus the distinction also occurs for
free minimal Cantor $\mathbb Z$-actions. The complete-corner
projection in a minimal realization belongs to the completed
diagonal and need not be continuous on the original Cantor space.

\subsection{The range for full invariant families of minimal actions}\label{subsec:boolean-actions}

\begin{thm}[Boolean-core action criterion]\label{thm:boolean-action-range}
Let $(\B,c)$ satisfy Theorem~\ref{thm:boolean-range}. It is isomorphic,
with capacity, to the completed projection algebra of a free
minimal Cantor $\mathbb Z$-action with its full invariant family
if and only if some countable $d_c$-dense Boolean subalgebra
$\B_0$ satisfies:
\begin{enumerate}
\item For $0\ne a\in\B_0$, $\nu(a)>0$ for all $\nu\in K_c$.
\item For every $\eps>0$, there is a finite partition of unity
$a_1,\ldots,a_n\in\B_0$ with $c(a_i)<\eps$ for all $i$.
\item If $a,b\in\B_0$ satisfy $\nu(a)<\nu(b)$ for every
$\nu\in K_c$, some $a'\in\B_0$, $a'\leq b$, satisfies
$\nu(a')=\nu(a)$ for all $\nu\in K_c$.
\end{enumerate}
\end{thm}
\begin{proof}
Restrict $K_c$ to probabilities on $Y=\Stone(\B_0)$.
It is compact convex. Condition (i) gives full support and
(ii) gives atomlessness of every measure.
These conditions also exclude Boolean atoms in $\B_0$:
any atom has positive mass under a fixed measure, whereas
every partition must contain it in one piece. Thus $Y$ is
Cantor. Condition (iii) is the Glasner--Weiss condition.
The dynamical-simplex theorem of Ibarluc\'ia--Melleray, with
Melleray's removal of the additional approximate-divisibility
hypothesis, supplies a minimal homeomorphism with exactly
this family \cite{IM,Melleray}. Its completion has $(\B,c)$
by Theorem~\ref{thm:boolean-range}; minimality on an infinite Cantor
space implies freeness.

Conversely, take the clopen core of the action.
Theorem~\ref{thm:boolean-recovery} identifies its full invariant family
with $K_c$. Minimality gives full support.
For each invariant measure choose a finite clopen partition
with all masses below $\eps$. The corresponding neighborhoods
cover the compact measure family; finitely many suffice, and
their common refinement gives (ii). The Glasner--Weiss
property gives (iii) \cite{IM,Melleray}.
\end{proof}

Full support is required on the core only, since completed
projections may vanish at some designated traces.
The three conditions imply that the recovered family is a
Choquet simplex; no separate simplex hypothesis is needed.

The bare Boolean algebra reconstructs the unmarked norm algebra
of an RR0 diagonal by Proposition~\ref{prop:boolean-Stone}.
The exact capacity additionally recovers the designated trace
simplex in the full-trace AF class. Bare classification is coarser:
all finite-vertex simplexes give the standard atomless probability
algebra, whereas Proposition~\ref{prop:boolean-different}
distinguishes two infinite-vertex examples. A criterion deciding
bare Boolean isomorphism for arbitrary different infinite
simplexes remains open in this discussion.

\end{document}